\input amstex
\documentstyle{amsppt}

\magnification=\magstep1 
\pagewidth{6.5truein}
\pageheight{9.0truein}

\document
\def\Bbb{\fam\msbfam \tenmsb}
\def\P{{\Cal P}}
\def\la{{\langle }}
\def\ra{{\rangle }} 
\def\k {{\bold k}}
\def\Om{{\Omega}}
\def\K/t{{K_n / \la T \ra}}
\def\d,{{\dots ,}}
\def\spec{\operatorname{spec}}
\def\Hspec{H\operatorname{-spec}}
\def\blt{{\bullet}}
\def\got{{\gamma_{1,2}}}

\def\Aut{\operatorname {Aut}}
\def\prim{\operatorname {prim}}
\def\height{\operatorname {ht}}
\def\GKdim{\operatorname {GK.dim}}
\def\AD{{\bf 1}}
\def\BG{{\bf 2}}
\def\FRT{{\bf 3}}
\def\GK{{\bf 4}}
\def\GKB{{\bf 5}}
\def\G{{\bf 6}}
\def\GLen{{\bf 7}}
\def\GLet{{\bf 8}}
\def\GW{{\bf 9}}
\def\HK{{\bf 10}}
\def\JZ{{\bf 11}}
\def\Mal{{\bf 12}}
\def\Man{{\bf 13}}
\def\YMan{{\bf 14}}
\def\MR{{\bf 15}}
\def\Mus{{\bf 16}}
\def\S{{\bf 17}}
\def\Q{{\bf 18}}
\def\OhP{{\bf 19}}
\def\Res{{\bf 20}}
\def\RTF{{\bf 21}}

\topmatter
\title The Prime and Primitive Spectra of Multiparameter Quantum
Symplectic and Euclidean Spaces 
\endtitle
\rightheadtext {Quantum Symplectic and Euclidean Spaces}

\author {K. L. Horton} \endauthor
\address Department of Mathematics, University of California, Santa
Barbara, CA  93106, USA
\email horton$\@$math.ucsb.edu \endemail
\endaddress

\abstract We investigate a class of algebras that provides multiparameter
versions of both quantum symplectic space and quantum Euclidean
$2n$-space.  These algebras encompass the graded quantized Weyl algebras,
the quantized Heisenberg space, and a class of algebras introduced by
Oh.  We describe the structure of the prime and primitive ideals of
these algebras.  Other structural results include normal separation and
catenarity.  
\endabstract

\endtopmatter

\bigskip

\head 
{Introduction}  
\endhead

The quantized coordinate rings known for affine spaces and quantum 
matrices have already been introduced in multiparameter versions, as
studied in \cite{\Man}, \cite{\YMan}, and \cite{\Res}, for instance. 
While single parameter quantum symplectic and Euclidean spaces have been
studied in \cite{\Mus} and \cite{\RTF},  no multiparameter versions of
these quantum spaces have been explicitly given.  We
have worked out a class of multiparameter algebras which is broad enough
to fit the pattern of the generators and the relations in both the
quantum symplectic and the quantum Euclidean $2n$-spaces.  The algebras
in this class are the most general algebras fitting this pattern that
are also iterated skew polynomial rings \cite {\HK}.  Further, they
incorporate the graded quantized Weyl algebras, the quantum Heisenberg
space \cite {\JZ}, and the algebras studied by Oh in \cite{\S}.

There exist natural tori that act as automorphims on these algebras, and
it is  known from the work of Brown, Goodearl, and Lenagan that the key
to understanding the prime and primitive ideals of the algebras is to
pin down the prime ideals invariant under these automorphisms.  We
completely determine these invariant prime ideals; this  generalizes the
earlier work of Oh in \cite {\S} and \cite {\Q}, as well as that of 
G\'omez-Torrecillas, El Kaoutit, and Benyakoub in
\cite {\GKB} and \cite {\GK}, using different methods.  In the final
section, we will illustrate how this leads to a complete determination
of all of the primitive ideals.  Other consequences include normal
separation and catenarity.  

Throughout, $\bold k$ will represent a base field of
arbitrary characteristic. For most results, $\bold k$ will need to
contain non-roots of unity.  All algebras will be unital.   
     
This work will form a portion of the author's PhD thesis at the University of California,
Santa Barbara.  

 \head
{1. The algebra $K_n = K_{n, \Gamma}^{P,Q}(\k).$} 
\endhead

In this section, we will define our algebras and show how they include
the quantum symplectic space and the other algebras that have been
previously studied.  We also provide some definitions and
observations that will be useful later.

\definition {Definition 1.1}  Let $\bold k$ be a field and let $P, Q \in (\k^{\times})^n$ such that
$P=(p_1, \dots p_n)$ and $Q = (q_1, \dots q_n)$ where $p_i \neq q_i$ for each $i \in \{1, \dots ,n\}.$  
Further, let $\Gamma = (\gamma_{i, j}) \in M_n(\bold k^{\times})$ with $\gamma_{j,i} = \gamma_{i,j}^{-1}$ 
and $\gamma_{i,i} =1$ for all $i,j.$   Then $K^{P,Q}_{n, \Gamma}(\bold k)$ is the $\bold k-$algebra generated 
by $x_1, y_1, \dots, x_n, y_n$  satisfying the following relations: 

$$\xalignat2 y_iy_j & = \gamma_{i,j}y_jy_i, &&\forall i,j \\
x_iy_j &= p_j\gamma_{j,i}y_jx_i &&(i<j) \\
x_iy_j &= q_j \gamma_{j,i}y_jx_i && (i>j) \\
x_ix_j &= q_ip_j^{-1}\gamma_{i,j}x_jx_i && (i<j) \\
x_iy_i &= q_iy_ix_i + \sum _{\ell < i}(q_{\ell}-p_{\ell}) y_{\ell}x_{\ell} &&\forall i.\endxalignat$$

\noindent When convenient, we will drop the $P,Q, \Gamma$ notation and write $K_n$ for $K^{P,Q}_{n, \Gamma}.$
\enddefinition

\definition {Definition 1.2}  Let 
$$H_n = \{ (h_1, h_2, \dots, h_{2n-1}, h_{2n}) \in (\k^{\times})^{2n} \mid h_{2i-1} h_{2i} = h_{2j-1} h_{2j} \hskip.05in 
\forall i, j = 1, \dots , n \}.$$

\noindent The group $H_n$ acts on $K_n$ by $\k$-automorphisms as follows:  for $h = (h_1, h_2, \dots, h_{2n-1}, h_{2n}),$ 
we have $h(x_i) = h_{2i-1} x_i$ and $h(y_i) = h_{2i}y_i$ for each $i \in \{1, \d, n \}.$   We will often drop the subscript 
and write $H$ for $H_n.$  
\enddefinition    
\medskip

Specific choices of $P,$ $Q,$ and $\Gamma$ will give rise to five algebras which have been previously studied. 

\definition {Example 1.3}  Given a nonzero element $q$ of $\k,$ we set $q_i = q^{-2}$ 
for each $i$ and $p_j = 1$ for each $j.$  Further, setting $\gamma_{i,j} =q$ whenever $i < j $ yields the 
{\it coordinate ring} ${\Cal {O}_q(\frak{sp}}(\k^{2n}))$ {\it of the quantum symplectic space\/}, 
the $\k$-algebra generated by $x_1, y_1, \dots, x_n, y_n$ with the following relations:

$$\xalignat2 y_iy_j & = q y_jy_i, &&(i<j)\\
x_ix_j & = q^{-1}x_jx_i, && (i < j)\\
x_iy_j & = q^{-1}y_jx_i, && (i \neq j)\\
x_iy_i & = q^{-2}y_ix_i + \sum_{\ell < i} (q^{-2} -1)y_{\ell}x_{\ell} && \forall i. \endxalignat$$

\smallskip

\noindent Faddeev, Reshetikhin, and Takhtadzhyan defined ${\Cal {O}_q(\frak{sp}}(\k^{2n}))$ in \cite{\FRT} 
and Musson gave new relations for the algebra in \cite{\Mus}.  Oh studied the primitive ideals of 
${\Cal {O}_q(\frak{sp}}(\k^{2n}))$ in \cite{\Q}, and the algebra is
considered in the latter two  papers as having generators $X_1,
X_{1^{\prime}}, \d, X_n, X_{n^{\prime}}.$  Setting $x_i =
X_{i^{\prime}}$  and $y_i = q^iX_i$ for each $i \in \{1, \d, n \}$
yields the relations given above.   G\'omez-Torrecillas, El Kaoutit, and
Benyakoub described a stratification of the spectra of 
${\Cal {O}_q(\frak{sp}}(\k^{2n}))$ via a torus of rank $n$ in \cite{\GKB}.  The generators in 
\cite{\GKB} are $Y_1, \d, Y_n, X_1, \d, X_n$ and are given from $K_n$ by setting $x_i = Y_i$ and $y_i =
q^iX_i.$

\enddefinition

\definition {Example 1.4}  Next, consider the case where $p_i= 1$ for each $i.$  Without further restriction on 
$\Gamma,$ the relations for $K_n$ become the relations for the 
{\it graded quantized Weyl algebra, \/}  gr$A_n^{Q, \Gamma}(\k):$

$$\xalignat2 y_iy_j &= \gamma_{i,j}y_jy_i, && \forall i, j\\
x_ix_j & = q_i \gamma_{i,j} x_jx_i, && (i < j)\\
x_iy_j & = \gamma_{j,i} y_jx_i, && (i < j)\\
x_iy_j & = q_j \gamma_{j,i}y_jx_i, && (i > j)\\
x_iy_i & = q_iy_ix_i + \sum_{\ell < i} (q_{\ell} -1)y_{\ell}x_{\ell} && \forall i. \endxalignat$$

\smallskip

\noindent The quantized Weyl algebra $A_n^{Q, \Gamma}(\k)$ arose in \cite{\Mal} and was further studied in 
\cite{\AD}.
\enddefinition

\definition {Example 1.5}  If $q \in \k^{\times},$ then setting $q_i = 1$ and $p_j = q^{-2}$ 
for each $i,j$ with $\gamma_{i,j} =q^{-1}$ for all $i < j$ forms the
{\it coordinate ring\/} $\Cal{O}_q(\frak{o} \k^{2n})$ {\it of quantum Euclidean 2n-space over\/} $\k:$  

$$\xalignat2 y_iy_j & = q^{-1}y_iy_j, && (i < j)\\
x_iy_j & = q^{-1} y_jx_i, && (i \neq j)\\
x_ix_j & = q x_jx_i, && (i < j)\\
x_iy_i & = y_ix_i + \sum_{\ell < i} (1 - q^{-2})y_{\ell}x_{\ell} && \forall i. \endxalignat$$

\smallskip
\noindent The algebra $\Cal{O}_q(\frak{o} \k^{N})$ arose in \cite{\FRT} and was given a simpler set of 
relations in \cite{\Mus}.  The generators for the even case with $N= 2n$ are given by
 $X_1, X_{1^{\prime}}, \d, X_n, X_{n^{\prime}}.$  
Setting $x_i = X_{(n+1-i)'}$ and $y_i = q^{n+1-i} X_{n+1 -i}$ for $i = 1, \d, n$ yields the above relations.  
Oh and Park studied the primitive ideals of $\Cal{O}_q(\frak{o} \k^{N})$ for both the even and odd cases in 
\cite{\OhP}. 
\enddefinition

\definition {Example 1.6}   Next, suppose that $q \in \k^{\times}.$  Setting 
$p_i = q^2$ and $q_i = 1$ with $\gamma_{i,j} = q$ whenever $i < j$ gives rise
to the {\it coordinate ring of quantum Heisenberg space \/} $F_q(n):$

$$\xalignat2 y_iy_j & = q y_jy_i, && (i < j)\\
x_iy_j & = q y_jx_i, && (i \neq j)\\
x_ix_j & = q^{-1} x_jx_i, && (i < j)\\
x_iy_i & = y_ix_i + \sum_{\ell < i} (1 - q^2)y_{\ell} x_{\ell}&&  \forall i. \endxalignat$$           

\noindent  The quantum Heisenberg space was first introduced by Faddeev, Rashetikhin, and Takhadjian in 
\cite{\FRT} and Jacobsen and Zhang studied $F_q(n)$ in \cite{\JZ} where $\k = \Bbb {C}$ and $q$ is a root of unity. 
In that paper, they considered generators $z_0, \dots , z_{n-1}, z_1^* \dots , z_{n-1}^*$ for $F_q(n).$
Setting $z_i = y_{n-i}$ and $z_i^* = x_{n-i}$ for each $i$ gives the equivalent algebra above.   
\enddefinition
\medskip

\noindent $\bold {Example}$ $\bold {1.7.}$  Lastly, let  $\lambda, d  \in \k^{\times}$ and set $p_i = d^{-1}$ for each
$i$ without further restrictions on $Q$ and $\Gamma.$  Then the relations for  $K_n$ become:

$$\xalignat2 y_iy_j &= \gamma_{i,j} y_jy_i, && \forall i, j\\
x_iy_j &= d^{-1} \gamma_{j,i} y_jx_i && (i<j)\\
x_iy_j &= q_j \gamma_{j,i} y_jx_i && (i>j)\\
x_ix_j &= q_i d \gamma_{i,j} x_jx_i && (i < j)\\
x_iy_i &= q_i y_ix_i + \sum_{\ell < i} (q_{\ell} - d^{-1})y_{\ell}x_{\ell} && \forall i. \endxalignat$$

\noindent Set $C= (Q,d, \lambda, 0)$ and define $X_i = d^ix_i$ and $Y_i = \lambda^i y_i$ for each $i.$ 
Then $K_n$ gives the $R^{C, \Gamma} _n (\k)$ algebra introduced by Oh in
\cite{\S}.  This algebra was further  studied by G\'omez-Torrecillas
and El Kaoutit, who classified its prime and primitive ideals in
\cite{\GK}.   Note that the additional coefficients $p_1, \d, p_n$ allow
the $K_n$ algebra to cover more cases than than  the $R^{C,
\Gamma}_n(\k)$ examples for $C \in
((\k^{\times})^{n+1},\k^{\times},\k^{\times},0).$

\medskip
We will now consider a group $H$ acting by automorphisms on rings or $\k$-algebras.  If $H$ acts on two
rings $A$ and $B,$ a map $\phi: A \rightarrow B$ is said to be {\it $H$-equivariant\/}
if and only if $\phi (h(a)) = h (\phi (a))$ for each $h \in H,$ $a \in A.$  When
$\phi$ is an isomorphism, we write $A \cong_H B.$  An {\it $H$-eigenvector $x$ \/} of a $\k$-algebra
$A$ is a nonzero element $x \in A$ such that $h(x) \in \k^{\times} x$ for
each $h \in H.$   Note that the generators $x_1, y_1, \d, x_n, y_n $ are
$H_n$-eigenvectors of $K_n.$  

Whenever $H$ acts on a ring $R,$ an ideal $Q$ of $R$ is said to be $H$-{\it stable \/} if
$h(Q)=Q$ for all $h \in H.$   Further, a proper ideal $Q$ of $R$ is
$H$-{\it prime\/} if $Q$ is $H$-stable such that whenever $I,J$ are
$H$-stable  ideals of $R$ with $IJ \subseteq Q,$ either $I \subseteq Q$
or $J \subseteq Q.$  As in the usual case,      an $H$-{\it prime ring
\/} is a ring in which $0$ is an $H$-prime ideal.  A ring $R$ is said to
be
$H$-{\it simple\/} if $0$ and $R$ are the only $H$-stable ideals of $R.$

Recall that $z \in R$ is {\it normal\/} if $zR = Rz.$  Further, $r \in R$ is said to {\it normalize \/}
a subring $S$ of $R$ if $rS = Sr.$  Note that if $z$ is normal in $R,$ then $zR = Rz = \la z \ra.$

When considering skew polynomial rings, we will utilize left-hand coefficients.  That is, given 
a skew polynomial ring $S = R[x; \sigma, \delta],$ we have $xr = \sigma(r) x + \delta (r)$ for $r \in R$
instead of $rx = x \sigma(r) + \delta(r)$.  Further, $\sigma$ will always represent an automorphism.  
As observed in \cite{\G}, if $I$ is an ideal of $S$ such that $\sigma(I)
\subseteq I$ and $\delta(I) \subseteq I,$  then $(\sigma, \delta)$
induces a skew derivation on $R/I$ and $IS = SI$ is an ideal of $S$ with 
$S/IS \cong (R/I)[x; \sigma, \delta].$  

Any nonzero $s \in S$ may be written uniquely as 
$s = r_mx^m + r_{m-1}x^{m-1} + \dots r_ + r_1x + r_0$ for some $m \in {\Bbb Z}^+, r_i \in R$ with $r_m \neq 0.$  
Here, $m$ is the {\it degree \/} of $s$ and $r_m$ is the {\it leading coefficient\/} of $s.$  When 
$R$ is a domain, $R[x; \sigma, \delta]$ is a domain and deg$(sw) =$deg$(s) +$deg$(w)$ for each nonzero 
$s,w \in R[x; \sigma, \delta].$  Further, if $A = B[x; \sigma][y; \tau]$ for some $\k$-algebra $B$ with 
$\tau (B) = B$ and $\tau(x) \in \k ^{\times} x,$ then 
$A \cong B[y; \tau ^{\prime}][x; \sigma ^{\prime}]$ for some $k$-automorphisms
$\tau ^{\prime}, \sigma^{\prime}.$  More properties of skew polyomial rings can be found in \cite{\G}.
   
Another basic fact that will be useful is the following:  if $R$ is a $\k$-algebra with $X \subseteq R$ a 
multiplicative set, then $X$ is a right denominator set if and only if $\k^{\times} X$ is a right denominator set.

Many of the properties of prime ideals carry over to the $H$-prime case.  In particular, we will make use of two 
observations, with details left to the reader.
\smallskip

\definition {Observation 1.8}  If $Q$ is an $H$-prime ideal of a $\k$-algebra $A,$ then whenever $x$ and 
$y$ are $H$-eigenvectors in $A$, we have the following:

(i)  $xAy \subseteq Q \Rightarrow$ $x \in Q$ or $y \in Q$ and 

(ii) if either $x$ or $y$ is normal modulo $Q,$ then $xy \in Q$ implies that $x \in Q$ or $y \in Q.$         
\enddefinition

\definition {Observation 1.9}  Let $X$ be a right denominator set in a right noetherian ring $R,$
and suppose that $H$ acts on $R$ such that $X$ is $H$-stable.  Then the action of $H$ on $R$ extends uniquely 
to an action on the localization $RX^{-1}$ by automorphisms.  Further, extension and contraction provide inverse
bijections between the set of $H$-prime ideals of $RX^{-1}$ and the set of those $H$-prime ideals of $R$ that are
disjoint from $X.$  
\enddefinition

\definition {Remark 1.10}  The results of Observation 1.9 also follow if $R$ is a $\k$-algebra, $H$ acts 
by $\k$-algebra automorphisms, and we only assume that $\k^{\times}X$ is $H$-stable instead of $X.$  \enddefinition

\head
{2.  Admissible Sets}
\endhead

We will define admissible sets and show that the ideals that they
generate are both prime and $H$-prime.  From Lemma $2.1$ on, we will
assume that $p_iq_i^{-1}$ is not a root of unity for each $i \in \{1,
\dots , n \}.$

\proclaim {Lemma 2.1}  For each $i= 1, \dots ,n,$ let $\Omega_i =
\sum_{\ell \leq i}(q_{\ell}-p_{\ell})y_{\ell}x_{\ell},$ and set $\Om_0 =
0.$

{\rm (a)}  For any $\Omega_i,$ 

$$\xalignat2 \Om_i x_j & = p_j^{-1} x_j \Om_i, && 1 \leq i < j \leq n\\ 
\Om_i x_j & = q_j^{-1} x_j \Om_i, && 1 \leq j \leq i \leq n\\
\Om_i y_j & = p_j y_j \Om_i, && 1 \leq i < j \leq n\\
\Om_i y_j & = q_j y_j \Om_i, && 1 \leq j \leq i < n\\
\Om_i \Om_j & = \Om_j \Om_i && 1 \leq j \leq n. \endxalignat $$

As a result, $\Om_i$ is normal in each of $K_i, \d, K_n$ and normalizes each of $K_1, \dots, K_{i-1}.$     

{\rm (b)}  We have the following relations:   
$$\xalignat1 \Omega_{i-1} & = x_iy_i - q_iy_ix_i\\
\Omega_i & = x_iy_i -p_iy_ix_i. \endxalignat$$    

Thus, the cosets of $x_i$ and $y_i$ are normal in the quotient algebras $K_n/ \la \Om_{i-1} \ra$ and
$K_n/ \la \Om_i \ra.$
\endproclaim  

\smallskip

\noindent $\bold {Proof:}$  The above formulas may be derived 
from the definition of $K_n.$  For the last part of (b), recall that
$x_i$ and $x_j$ commute up to scalar multiplication for all $i, j\in
\{1,\dots ,n\}$, and that the same is true for $y_i$ and $y_j.$ Since
$x_i$ and $y_{\ell}$ also commute up to scalar multiplication for $i
\neq \ell$ and $\overline {x_i y_i} =
 q_i \overline {y_i x_i}$ with $\overline {y_i x_i} = q_i^{-1}\overline
{x_iy_i} $ in $K_n / \la \Om_{i-1} \ra,$ we have that $\overline x_i$ and
$\overline y_i$ are normal in $K_n / \la \Om_{i-1} \ra.$  The case for
$K_n / \la \Om_i \ra$ is similar.  $\blacksquare$   

\medskip

\definition {Definition 2.2}  Set $\P_n = \{x_1, \d, x_n, y_1, \d, y_n, \Omega_1, \d, 
\Omega_n\}.$  Following Oh in \cite {\Q}, a subset $T \subseteq \P_n$ is
{\it admissible\/} if $T$ satisfies the following two conditions:

(1) $x_i \in T$ or $y_i \in T$ if and only if $\Omega_i \in T$ and $\Omega_{i-1} \in T, \hskip.25in 2 \leq i \leq n.$

(2) $x_1 \in T$ or $y_1 \in T$ if and only if $\Omega_1 \in T.$ 

\enddefinition

\definition {Definition 2.3}  Given an admissible set $T,$ let $N_T$ be
the subset of $\P_n$  defined by the following conditions:

{\rm (a)}  $\Om_i \in N_T$ if and only if $  \Om_i \notin T;$

{\rm (b)}  $x_1 \in N_T$ if and only if $ x_1 \notin T;$

{\rm (c)}  $y_1 \in N_T$ if and only if $y_1 \notin T;$

{\rm (d)}  for $i > 1,$ $x_i \in N_T$ if and only if $ x_i \notin T$ 
and $\Om_{i-1} \in T$ or $\Om_i \in T;$ and

{\rm (e)}  for $i > 1,$ $y_i \in N_T$ if and only if $ y_i \notin T$ and
$\Om_{i-1} \in T$ or $\Om_i \in T.$
\enddefinition

\proclaim{Lemma 2.4} Let $P$ be an $H$-prime ideal of $K_n,$ and set $T = P \cap \P_n.$ 
Then $T$ is an admissible set.
\endproclaim
\smallskip

\noindent $\bold {Proof:}$  Suppose that $x_i \in T,$ where $1 \leq i \leq n.$  Then
$x_i \in P,$ so that $x_iy_i,$ $q_iy_ix_i,$ and $p_iy_ix_i \in P.$  For $i = 1,$
$\Omega_1 = (q_1 - p_1)y_1x_1 \in P.$  If $i > 1,$ $\Omega_{i-1} = x_iy_i - q_iy_ix_i
\in P,$ and $\Omega_i = x_iy_i -p_iy_ix_i \in P.$  Hence, for $i=1,$ if $x_i \in T,$
then $\Omega_i \in T,$  and for $i > 1,$ if $x_i \in T,$ then $\Omega_i, \Omega_{i-1} \in
T.$ Analogously, for any $i$ such that $1 \leq i \leq n,$ if  $y_i \in P,$ then
$\Omega_i \in T,$ and for $i >1,$  we have that $\Omega_{i-1} \in T$ whenever $y_i
\in T.$

Next, suppose that $\Omega_1 \in T.$  Then $(q_1 - p_1)y_1x_1 \in P,$ so
that $y_1x_1 \in P.$  Now, $x_1$ and $y_1$ are normal in $K_n,$ so
$y_1K_nx_1 = y_1x_1K_n \subseteq P,$ and by Observation 1.8, since $x_1$
and $y_1$ are $H$-eigenvectors, either $x_1 \in P$ or $y_1 \in P.$  That
is, if $\Omega_1 \in T,$ then either $x_1 \in T$ or $y_1 \in T.$
 
If $\Omega_i, \Omega_{i-1} \in T$ for some $i > 1,$ then each of  
$x_iy_i -p_iy_ix_i$ and $ x_iy_i - q_iy_ix_i$ is contained in $P,$ so
that $(p_i - q_i)y_ix_i \in P.$  Since $p_i - q_i \neq 0,$ it follows
that $y_ix_i \in P.$  Now, $\overline x_i$ and $\overline y_i$ are
normal in $\overline K_n = K_n/P$ since $\Omega_i \in P,$ so 
$\overline x_i \overline K_n \overline y_i =  \overline K_n \overline
x_i \overline y_i = 0.$ Further,  $\overline K_n$ is an $H$-prime ring
because $P$ an $H$-prime ideal.  By Observation 1.8, either
$\overline x_i= 0 $ or $\overline y_i = 0;$ that is, either $x_i \in P$
or $y_i \in P.$ $\blacksquare$

\smallskip

\proclaim{Proposition 2.5}  The algebra $K_n$ is an iterated skew polynomial ring.  Hence, $K_n$ is 
noetherian and an integral domain.  Further, there is a $\k$-basis for $K_n$ consisting of

$$ \Cal {A} = \{x_1^{r_1} y_1^{r_2} \dots x_n^{r_{2n-1}} y_n^{r_{2n}} \mid r_i \in {\Bbb {Z}}^+\}.$$
\endproclaim 

\noindent $\bold {Proof:}$  Note that $K_n$ is an iterated skew polynomial ring
$$K_n = \k[x_1] [y_1; \tau_1] [x_2; \sigma_2] [y_2; \tau_2, \delta_2] \dots [x_n; \sigma_n]
[y_n; \tau_n, \delta_n]$$
\noindent for automorphisms $\sigma_i, \tau_i$ and $\tau_i$-derivation $\delta_i$ defined as follows:

$\sigma_i: \k [x_1, y_1, \dots , x_{i-1}, y_{i-1}] \longrightarrow  \k
[x_1, y_1, \dots , x_{i-1}, y_{i-1}],$
 $$\xalignat2 \sigma_i (x_j) &= q_j^{-1} p_i \gamma_{i,j} x_j &&1 \leq j
\leq i -1,\\ \sigma_i (y_j) & = q_j \gamma_{j,i} y_j && 1 \leq j \leq
i-1,\endxalignat$$  

$\tau_i: \k [x_1, y_1, \dots , x_{i-1}, y_{i-1}, x_i] \longrightarrow \k
[x_1, y_1, \dots , x_{i-1}, y_{i-1}, x_i],$
$$\xalignat2 \tau_i(x_j) & = p_i^{-1} \gamma_{j,i} x_j && 1 \leq j \leq i-1,\\
\tau_i(y_j) & = \gamma_{i,j} y_j && 1 \leq j \leq i-1,\\
\tau_i(x_i) & = q_i^{-1} x_i, &&  \endxalignat$$

$\delta_i : \k [x_1, y_1, \dots , x_{i-1}, y_{i-1}, x_i] \longrightarrow \k [x_1, y_1, \dots , x_{i-1}, y_{i-1}, x_i],$
$$\xalignat2  \delta_i(x_j) & = 0 && 1 \leq j \leq i-1,\\
\delta_i(y_j) & = 0 && 1 \leq j \leq i-1,\\
\delta_i(x_i) & = -q_i^{-1} \sum_{\ell < i} (q_{\ell} - p_{\ell}) y_{\ell} x_{\ell}.  &&  \endxalignat $$

\noindent The remaining conclusions follow by standard results.  $\blacksquare$ 
\smallskip

\smallskip
\proclaim{Lemma 2.6}  Let $A = B[x; \sigma ] [y; \tau, \delta ],$ where $B$ is a domain
and a $\bold k$-algebra, and $\sigma, \tau$ are 
\noindent $\bold k$-automorphisms such that
$\tau (x) = \alpha x$ for some $\alpha \in \bold k^\times.$  Suppose further that
there exists an element of the form $\Omega = xy + z,$ with $z \in B$
such that $\Omega$ normalizes $B[x; \sigma]$ and is normal in $A$.  Then:
\smallskip

 {\rm (i)} if $z =0,$ then $\delta = 0.$

 {\rm (ii)} $A[x^{-1}]= B[x^{\pm1};\sigma] [x^{-1}\Omega; \tau]$.

{\rm (iii)} if $z \neq 0,$
 
    \hskip.3in {\rm (a)} $A[x^{-1}]\Omega \cap A= A\Omega$.

    \hskip.3in {\rm (b)} $A/\langle \Omega \rangle$ is a domain. 
\endproclaim

\noindent $\bold {Proof:} $ (i) If $z = 0,$ then $\Omega =xy$
normalizes $B[x; \sigma].$  For any $b \in B[x; \sigma],$ we have that
$xyb  = x \tau (b)y + x \delta (b)$ so that $x \tau (b)y +  x \delta (b)
\in xy B[x; \sigma] = B[x; \sigma]xy.$ Then $x \tau(b)y + x \delta(b) =
\hat bxy$ for some $\hat b \in B[x; \sigma]$ and hence, $x
\delta(b)=({\hat b}x - x \tau(b))y.$ As a result, 
$x \delta (b) = 0,$ and since $B[x; \sigma]$ is a 
domain, $\delta (b) = 0.$  It follows that $\delta = 0.$       

(ii) In the case where $z = 0,$ note that $A = B[x; \sigma][y; \tau].$ 
Then 
$$\alignat 1 A[x^{-1}] &= (B[x; \sigma][y; \tau])[x^{-1}] = B[x^{\pm 1}
\sigma][y; \tau]\\ &= B[x^{\pm 1}; \sigma][x^{-1}(xy); \tau] = 
B[x^{\pm1}; \sigma][x^{-1} \Omega; \tau]. \endalignat $$  

For $z \neq 0,$ it is well known that  $S = \{\beta x^{i} \mid \beta
\in \bold k ^{\times}, i \in \Bbb {Z} ^+\}$ is a denominator set in
$B[x; \sigma],$ and that $B[x; \sigma][S^{-1}] = B[x^{\pm1};\sigma].$ 
Since $\tau(S) = S,$ we have that $S$ is also a denominator set in $A$ by
\cite{\G , Lemma 1.4}, which further yields that $A[S^{-1}]= B[x;\sigma]
[S^{-1}][y; \tau,\delta] =  B[x^{\pm1}; \sigma][y;\tau,
\delta].$  Observe that the set $X = \{x^i \mid i \in \Bbb {Z} ^+ \}$ is
a denominator set of both $B[x; \sigma]$ and $A.$  Further, $A[S^{-1}] =
A[X^{-1}] = A[x^{-1}],$  since $A$ is a $\bold k$-algebra, so $A[x^{-1}]
=  B[x^{\pm1}; \sigma][y;\tau, \delta].$  Unless otherwise noted, we will
consider the degree of an  element of A to be its degree as a polynomial
in $y.$ 

Since $\Omega$ normalizes $B[x; \sigma],$ we have $\Omega x = t \Omega$
for some $t \in B[x; \sigma].$  That is, $(xy + z)x = t(xy + z),$ so
$\alpha x^2y + x \delta (x) + zx = txy + tz.$  Comparing the leading
coefficients and cancelling an  $x$ yields that $\alpha x = t.$ 
Therefore, $\Omega x = \alpha x \Omega.$

Consider $x^{-1} \Omega = y + x^{-1}z,$ where $x^{-1}z \in B[x^{\pm1};
\sigma].$  Since $\Omega$ normalizes $B[x; \sigma]$ and $\Omega x =
\alpha x \Omega,$ we have that $ \Omega$ normalizes $B[x^{\pm1};
\sigma].$  Further, $x^{-1}$ normalizes $B[x^{\pm1};
\sigma],$ so $(x^{-1} \Omega ) B[x^{\pm1}; \sigma ]$ $ =  B[x^{\pm1};
\sigma ](x^{-1} \Omega).$  Consequently, there exists a $\k$-algebra
automorphism $\psi: B[x^{\pm1}; \sigma] \rightarrow B[x^{\pm1}; \sigma]$
such that $(x^{-1} \Omega)f = \psi (f) (x^{-1} \Omega)$ for all $f \in
B[x^{\pm1}; \sigma].$   Now, $(x^{-1} \Omega )f =  yf + x^{-1} z f =
\tau (f) y + \delta (f)+ x^{-1} z f$ and $\psi (f) (x^{-1} \Omega ) =
\psi (f) y + \psi (f) x^{-1}z.$  Comparing leading coefficients, we have
that $\psi(f) = \tau (f);$  It follows that $\psi = \tau.$   Noting that
$x^{-1}\Omega $ is monic of degree one, we have that its powers form a
basis for $A[x^{-1}]$ as a free left $B[x^{\pm1}; \sigma ]$-module. 
Hence, $A[x^{-1}] = B[x^{\pm1}; \sigma][x^{-1} \Omega;
\tau].$

(iii)  (a) We will first show that $x$ is not a left zero-divisor
modulo $\langle \Omega \rangle;$ that is, for any $f \in A,$ if $xf \in
\langle \Omega \rangle,$ then $f \in \langle \Omega \rangle.$  Note that
$A = B[x; \sigma][y; \tau, \delta]$ is a domain, and that deg$(fg)$ =
deg$(f)$ + deg$(g)$ for all $f, g \in A.$  Moreover, $\langle \Omega
\rangle = \Omega A$ since $ \Omega $ is normal in $A,$ so deg$(g)$ $\geq
1$ for all nonzero $g \in  \langle \Omega \rangle.$  If $f \in A$ such
that deg$(f) = 0,$ then deg$(xf) = 0.$  Thus, if $xf \in 
\langle \Omega \rangle$ we have that $ xf = 0,$ so $f = 0 \in \langle
\Omega \rangle.$    

Suppose that $m >0$ such that for each $h \in A,$ if deg$(h) < m$ and
$xh \in \langle \Omega \rangle,$ then $ h \in  \langle \Omega
\rangle.$  Let $f \in A$ such that deg$(f) = m$ and $xf \in  \langle
\Omega \rangle.$  Then $f = f_0 + f_1y + \dots +f_my^m$ for some $f_i
\in B[x; \sigma], f_m \neq 0,$ and since $\Omega$ is normal, there
exists some $g = g_0 + g_1y + \dots +g_sy^s,$ for $g_i \in B[x; \sigma],
g_s \neq0,$ such that $xf = \Omega g.$  That is, 
$$ \align xf_0 + xf_1y + \dots +xf_my^m & = (xy +z)g_0 + (xy + z)g_1y +
\dots  + (xy + z)g_sy^s \\
 & = x\tau(g_0)y + x\delta(g_0) + zg_0+ x\tau(g_1)y^2 + x \delta(g_1)y +
zg_1y \\
& \qquad \qquad \qquad+ \dots +x \tau(g_s)y^{s+1} + x\delta (g_s)y^s +
zg_sy^s\\ 
& = x\delta(g_0) + z g_0 + [x\tau(g_0) + x \delta (g_1) + zg_1]y + \dots
+ x \tau (g_s)y^{s+1}. \endalign$$
   
Comparing the leading coefficients, we have that $s + 1 = m,$ so that
deg$(g') = m -2,$ where $g' = g_1 + g_2y + \dots +g_sy^{s-1},$ with $g =
g_0 + g'y.$  Considering the terms of degree zero in the equation $xf =
\Omega g,$ we have that $xf_0 = x \delta (g_0) + z g_0,$ or $x(f_0 -
\delta(g_0)) = zg_0.$  Since each of $f_0,$ $\delta(g_0),$ $z,$ and $g_0$
is an element of $B[x; \sigma],$ we may consider $x(f_0 - \delta(g_0))$
and $zg_0$ as polynomials in $x.$  Then $x(f_0 - \delta(g_0))$ has zero
constant term (we are including the possibility that $f_0 - \delta(g_0)
= 0$), and thus, $zg_0$ has zero constant term. Since $z$ is a regular
element of $B,$ it follows that $g_0 = xh_0$ for some $h_0 \in B[x;
\sigma],$ and we may write $g = xh_0 + g'y.$  Hence, 
$xf = \Omega(xh_0 + g'y) = \Omega x h_0 + \Omega g'y = \alpha x
\Omega h_0 + \Omega g'y,$ and $x(f - \alpha \Omega h_0) = \Omega g'y.$        

Observe that $\Omega g'y$ has zero constant term as a polynomial in
$y;$ then $f -\alpha \Omega h_0$ has zero constant term, and may be
written as $f - \alpha \Omega h_0 = f'y$ for some $f' \in A.$  Further,
$xf'y = \Omega g'y$ implies that $xf' = \Omega g'.$  Since deg$(g') = 
m -2,$ it follows that deg$(f') = m - 1.$  By our induction hypothesis,
$f'\in \langle \Omega \rangle,$ and hence, $f - \alpha \Omega h_0 \in
\langle \Omega \rangle,$ so $f \in \langle \Omega \rangle.$ 
This completes the induction step, yielding that $x$ is not a left
zero-divisor modulo $\langle \Omega \rangle. $ 
 
Consider $I = A[x^{-1}] \Omega;$ if $g \in I \cap A,$ there exists $h
\in A[x^{-1}]$ such that $g = h \Omega$ and, letting $m$ be the largest
power of $x^{-1}$ in $h,$ we have that $x^mg = x^mh\Omega,$ where $x^mh
\in A.$  That is, $x^mg \in \langle \Omega \rangle,$ the ideal generated
by $\Omega$ in $A.$  Since $x$ is not a left zero-divisor modulo
$\langle \Omega \rangle,$ it follows that $g \in \langle \Omega
\rangle.$  Hence, $I \cap A \subseteq \langle \Omega \rangle$ so that
$I \cap A = \langle \Omega \rangle = A \Omega.$
 
(b) Lastly, note that $I = A[x^{-1}] \Omega = A[x^{-1}](x^{-1} 
\Omega)$ since $x^{-1}$ is invertible.  As in the proof of (a) above,
$x^{-1} \Omega$ is normal in $A[x^{-1}],$ so $I =(x^{-1} \Omega)
A[x^{-1}],$ or $I =\langle x^{-1} \Omega \rangle$ in $A[x^{-1}].$  Then
$A[x^{-1}]/I = A[x^{-1}]/ \langle x^{-1} \Omega \rangle \cong
B[x^{\pm1}; \sigma],$ a skew polynomial ring over a domain.  By (a), $A/
\langle \Omega \rangle$ embeds in $A[x^{-1}]/I,$ and is thus a domain.
$\blacksquare$ 
\smallskip

\proclaim{Lemma 2.7 }  Let $A= B[x; \sigma][y; \tau, \delta],$ where $B$ 
is a domain and a $\k$-algebra, and let $\sigma \text{ and } \tau$ be
$\k$-automorphisms such that $\tau(x) = \alpha x$ for some 
$\alpha \in \k^{\times}.$  Suppose further that $\delta (B) = 0$ and $\delta(x) \in B$ such that
$\delta(x)$ is normal in both $A$ and $B.$  If $\delta (x) $ is an eigenvector of both 
$\sigma$ and $\tau,$ and if the quotient algebra $B/ \delta(x)B$ is
nonzero, then:

{\rm (i)} \hskip.15in $x \notin \la y \ra .$

{\rm (ii)} \hskip.1in $ y \notin \la x \ra.$

\endproclaim  

$\bold {Proof:}$  (i)  Set $z = \delta (x)$ and $C = B[x; \sigma].$ 
Then $zC = Cz,$ so $zC$ is the ideal of $C$ generated by $z.$  Similarly,
$Bz = zB$ and $zA = Az$  are the ideals  generated by $z$ in $B$ and
$A,$ respectively.  Since $z$ is an eigenvector of $\sigma,$ we have
that $\sigma (Bz) = Bz.$  As observed in \cite{\G}, $C/ zC \cong
(B/zB)[\overline x; \overline \sigma].$ 

Next, note that $\tau (zC) = zC$ and $\delta (z) = 0,$ so $\delta (zC)
\subseteq zC.$ Applying the observation  a second time yields that 
$$A/ zA \cong (C/zC)[\overline y; \overline \tau, \overline \delta]
\cong (B/zB) [\overline x; \overline \sigma][\overline y; \overline
\tau, \overline \delta].$$   Since $\delta (B) = 0$ and $\delta (x) \in
zB,$ we have that $\overline \delta = 0$ and hence,
$$A/zA \cong (B/zB)[\overline x; \overline \sigma][\overline y;
\overline \tau].$$  If $x \in \la y \ra,$ the ideal of $A$ generated by
$y,$ then $\overline x = c \overline y$ for some $c \in A/zA,$ a
contradiction of the skew polynomial ring  construction.  Thus, $x
\notin \la y \ra.$ 

(ii)  Considering $A/zA \cong (B/zB)[\overline x; \overline
\sigma][\overline y; \overline \tau],$ we note that $\overline
\tau(\overline x) = \alpha \overline x.$  Then observe that
$$(B/zB)[\overline x; \overline \sigma][\overline y; \overline \tau] = 
(B/zB)[\overline y; \hat \tau][\overline x; \hat \sigma]$$ for 
$\k-$algebra automorphisms $\hat \tau$ and $\hat \sigma.$ 
Consequently, $\overline y$ is not a multiple of $\overline x$ in
$A/zA,$ so $y \notin \la x \ra.$ $\blacksquare$

\smallskip
For ease of notation, we will allow $x_n, y_n$ to represent their
cosets in the factor algebras below.  

\proclaim{Theorem 2.8}  Let  $T$ be an admissible set of $K_n.$ 
Then $\langle T \rangle \cap \P_n = T$ and $\langle T \rangle$ is
completely prime.  Consequently, $\la T \ra$ is an $H$-prime ideal.
\endproclaim
\smallskip

\noindent $\bold {Proof:}$  We will proceed by induction on $n.$   If
$n = 1,$ there are four admissible sets, namely:  $\emptyset, \{y_1,
\Omega_1 \}, \{x_1, \Omega_1 \},$ and $\{x_1, y_1, \Omega_1 \};$ which 
respectively generate the ideals $0, \la y_1 \ra, \la x_1 \ra, $ and
$\la x_1, y_1 \ra.$ Clearly, $0 \cap \P_1 = \emptyset$ and $\la x_1, y_1
\ra \cap \{x_1, y_1, \Omega_1 \} = \{ x_1, y_1, \Omega_1 \}$.  
 Further, $\Om_1$ is an element of both $ \la x_1 \ra$ and 
$ \la y_1 \ra$ since $\Om_1 = (q_1 -p_1)y_1x_1.$ Noting that $y_1K_1 =
K_1y_1,$ we have that $\la y_1 \ra = y_1K_1$ and for each $a \in \la y_1
\ra,$ deg$(a) \geq 1$ as a polynomial in $y_1.$ Thus, $x_1 \notin \la
y_1 \ra$ and $\la y_1 \ra \cap \P_1 = \{y_1, \Omega_1 \}.$

Since $\tau_1(x_1) \in \k^{\times}x_1,$ we have that $K_1 \cong \k[y_1;
\tau_1'] [x_1; \sigma_1']$ for some $\k$-automorphisms $\tau_1',
\sigma_1'.$  Then $x_1K_1 = K_1x_1$ implies that each $b \in \la x_1
\ra$ satisfies deg$(b) \geq 1$ as a polynomial in $x_1.$  Consequently,
$y_1 \notin \la x_1 \ra,$ and $\la x_1 \ra \cap \P_1 = \{x_1, \Omega_1
\}.$  To see that $\la T \ra$ is completely prime  for each admissible
set $T$ of $K_1,$ note that the quotient algebra $K_1/  \la T \ra$ is 
isomorphic to either $K_1, \k[x_1], \k[y_1],$ or $\k;$ each of which is
a domain.        

Suppose now that $n > 1$ and that for each $\ell < n,$ if $U \subseteq
K_{\ell}$ is an admissible set and if $I$ is the ideal of $K_{\ell}$
generated by $U,$ then $I \cap \P_{\ell} = U$ and $I$ is completely
prime. Throughout, we will use $\la \dots \ra$ to denote an ideal of
either $K_n$ or a factor of $K_n,$ and will give names to ideals  in
other rings.  Given an admissible set $T$ of $K_n,$ let $T_{n-1} = T
\cap \P_{n-1}$ and note that $T_{n-1}$  is admissible as a subset of
$K_{n-1}.$  Further, each $a \in T_{n-1}$ is an eigenvector of both
$\sigma_n$ and $\tau_n,$ with  $\delta_n(a)=0.$   Letting $I_{n-1}$ be
the ideal of $K_{n-1}$ generated by $T_{n-1},$ we have that
$\sigma_n(I_{n-1}) = I_{n-1}$ and $\tau_n(I_{n-1}) = I_{n-1}$ with
$\delta_n(I_{n-1})=0.$

If $J$ is the ideal of $K_{n-1}[x_n; \sigma_n]$ generated by $T_{n-1},$
then, as observed in \cite{\G}, $$(K_{n-1}/I_{n-1})[\overline x_n;
\overline \sigma_n] \cong (K_{n-1}[x_n; \sigma_n])/J,$$ where $\overline
\sigma_n$ is the $\k$-algebra automorphism of $K_{n-1}/I_{n-1}$ induced
by $\sigma_n.$   Further, if $I_n$ is  the ideal of $K_n$ generated by
$T_{n-1},$ then $$K_n/I_n \cong  ((K_{n-1}[x_n; \sigma_n])/J)[\overline
y_n; \overline \tau_n, \overline \delta_n] \cong
(K_{n-1}/I_{n-1})[\overline x_n; \overline \sigma_n][\overline y_n;
\overline \tau_n, \overline \delta_n]$$  where $\overline \tau_n,
\overline \delta_n$ are induced by $\tau_n,\delta_n,$ respectively.  Set
$ A = (K_{n-1}/I_{n-1})[\overline x_n; \overline \sigma_n][\overline y_n;
\overline \tau_n, \overline \delta_n];$  by our induction hypothesis,
$K_{n-1}/I_{n-1}$ is a domain, so $A$ is a domain.  Note further that if
$a \notin I_{n-1},$  so that $\overline a \neq 0$ in $A,$ then $a \notin
I_n.$  Applying the induction hypothesis a second time yields that 
$I_n \cap \P_{n-1} = I_{n-1} \cap \P_{n-1} = T_{n-1}.$

Let $S_n = \P_n \setminus \P_{n-1} = \{x_n, y_n, \Om_n \},$ and $S = T
\setminus T_{n-1}.$  To show that $\la T \ra \cap \P_n = T,$ we will
first show that $\la T \ra \cap S_n = S,$ and then that $\la T \ra \cap
\P_{n-1} = T_{n-1}.$  There are five possibilities for $S,$
namely:  $\emptyset, \{x_n, \Om_n \}, \{y_n, \Om_n \}, \{ \Om_n \},$ and
$\{x_n,y_n, \Om_n \}.$  

If $S = \emptyset,$ note that $\overline x_n, \overline y_n$ are
nonzero by the skew polynomial ring construction. Thus,
$\overline x_n, \overline y_n \notin \la T \ra /I_{n-1},$ so $x_n, y_n
\notin \la T \ra.$  If $\Om_n \in \la T \ra  = \la T_{n-1} \ra,$ then
$\overline x_n \overline y_n - p_n \overline y_n \overline x_n = 0$ in
$A.$  But $\overline x_n \overline y_n = q_n \overline y_n \overline x_n
+ \overline \Om_{n-1},$ where $\overline \Om_{n-1} \in K_{n-1}/I_{n-1},$
so$\overline x_n \overline y_n = q_np_n^{-1} \overline x_n \overline y_n
+ \overline \Om_{n-1}$ implies that $q_np_n^{-1} = 1,$ a contradiction
since $p_n \neq q_n.$  Thus, $\Om_n \notin \la T \ra,$ and we have that
$\la T \ra \cap S_n = \emptyset = S.$

If $S = \{x_n, \Om_n \},$ note that $\Om_n = x_ny_n - p_ny_nx_n \in \la
x_n \ra,$ and $\la T \ra /I_n \cong  \la \overline x_n \ra,$ the ideal of
$A$ generated by $\overline x_n.$  By Lemma 2.7, $\overline y_n \notin \la
\overline x_n \ra,$ so $y_n \notin \la T \ra.$  Thus, $\la T \ra \cap
S_n = \{x_n, \Om_n \}.$

If $S = \{y_n, \Om_n \},$ then $\Om_n \in \la y_n \ra,$ and 
$\la T \ra /I_n \cong \la \overline y_n \ra,$ the ideal of $A$ generated
by $\overline y_n.$   Again by Lemma 2.7, $\overline x_n \notin \la
\overline y_n \ra$ in $A,$ so $x_n \notin \la T \ra.$   Thus,
$\la T \ra \cap S_n  = \{y_n, \Om_n \}.$

If $S = \{\Om_n \},$ then $\la T \ra /I_n \cong \la \overline \Om_n
\ra,$ the ideal of $A$ generated by $\overline \Om_n.$  Since $\la \Om_n
\ra \subseteq \la y_n \ra$ and $x_n \notin \la y_n \ra,$ we have that
$x_n \notin \la \Om_n \ra.$  Similarly, $y_n \notin \la  \Om_n  \ra,$
yielding that $\la T \ra \cap S_n = \{ \Om_n \}.$ 

Lastly, if $S = \{x_n, y_n, \Om_n \} = S_n,$ it is clear that $S_n \cap
\la T \ra = S.$  Thus, for each possible set $S,$ we have that $S_n \cap
\la T \ra =S.$ 

To show that $\la T \ra \cap \P_{n-1} = T_{n-1},$ we observe that $\la
T \ra = \la T_{n-1} \ra + \la S \ra = I_n + \la S \ra.$  As noted above,
$I_n \cap \P_{n-1} = T_{n-1}.$  Thus $(\la T_{n-1} \ra + \la S \ra)/I_n
= \la S \ra /I_n,$ and it suffices to show that $\la \overline S \ra
\cap K_{n-1}/I_{n-1} = 0$ in $A.$  If $S = \emptyset,$ then $\la S \ra =
0$ and the result is clear.  For  $S = \{x_n, y_n, \Om_n \},$ we have
that  $\Om_{n-1} \in T_{n-1}$ since $T$ is admissible, so $\Om_{n-1} \in
I_{n-1}.$  It follows from Lemma 2.1 that $\overline x_n$ and $\overline
y_n$ are normal in $A.$  

Now, $\la \overline S \ra = \la \overline x_n, \overline y_n
\ra,$ so if $w \in \la \overline S \ra \cap K_{n-1}/I_{n-1},$ then
$w = \overline x_n a + b \overline y_n$ for some $a, b \in A.$ 
By the skew polynomial construction, there exist $a_0, \dots, a_m, b_0,
\dots, b_r \in (K_{n-1}/I_{n-1})[\overline x_n; \overline \sigma_n]$
such that $a = \Sigma_{i=0}^m a_i \overline  y_n^i$ and $b =
\Sigma_{j=0}^r b_j \overline y_n^j.$  Then $w = \overline x_na + b
\overline y_n = \Sigma_{i=0}^m \overline x_n a_i \overline y_n^i + 
\Sigma_{j=0}^r b_j \overline y_n^{j+1},$ so $w = \overline x_n a_0.$

By the construction of $(K_{n-1}/I_{n-1})[\overline x_n; \overline
\sigma_n],$ there exist $ c_0, \dots, c_s \in K_{n-1}/I_{n-1}$  such
that $a_0 =  \Sigma_{i=0}^s c_i \overline x_n^i.$  Then  $w = \overline
x_n a_0 = \Sigma_{i=0}^s \overline x_n c_i \overline x_n^i =
\Sigma_{i=0}^s \overline \sigma(c_i) \overline x_n ^{i + 1}.$  Hence, 
$w = 0,$ yielding that  $\la \overline x_n, \overline y_n \ra\cap K_{n-1}
= 0,$ as desired. 

If $S = \{x_n, \Om_n \}$ or $S = \{y_n, \Om_n \},$ then $\Om_{n-1} \in
T_{n-1}$ and $\la S \ra \cap K_{n-1} = 0$ by  the above argument.  For
$S = \{ \Om_n \},$ note that $\overline {\Om}_n$ is a normal element of
$A.$  Thus, $A \overline \Om_n = \overline \Om_nA = \la \overline \Om_n
\ra$, the ideal generated by $\overline {\Om}_n.$  Further,
deg$(\overline {\Om}_n) = 1$ as a polynomial  in $\overline y_n,$ so
$\deg (a) \geq 1$ for each $a \in \la \overline {\Om}_n \ra.$ 
Consequently, $\la \overline \Om_n \ra \cap  K_{n-1}/I_{n-1} = 0$ and
hence, $\la T \ra \cap \P_n= T$ for any admissible set $T.$

To see that $\la T \ra$ is completely prime for any admissible set $T$
of $K_n,$ recall that the set $T_{n-1} =  T \cap \P_{n-1}$ generates a
completely prime ideal $I_{n-1}$ of $K_{n-1}$ by the induction
hypothesis.   As noted above, each element of $T_{n-1}$ is an
eigenvector of both $\sigma_n$ and $\tau_n$ with 
$\delta_n (T_{n-1}) = 0,$ so $I_{n-1}K_n=K_nI_{n-1}$ is an ideal of
$K_n.$   We then have five possibilities:  

If $\Omega_n \notin T,$ then $\la T \ra = I_{n-1}K_n,$ and $K_n/
\langle T \rangle \cong (K_{n-1}/I_{n-1})[x_n;\overline \sigma_n][y_n;
\overline \tau_n, \overline \delta_n] = A$  as above.  Again, since
$I_{n-1}$ is a completely prime ideal of $K_{n-1},$  it follows that
$(K_{n-1}/I_{n-1})[x_n;\overline \sigma_n] [y_n; \overline \tau_n,
\overline \delta_n]$ is a domain as an iterated skew polynomial ring
over a domain.  Hence, $K_n/\langle T \rangle$ is a domain.

If $x_n \in T$ but $y_n \notin T,$ then $K_n/ \langle T \rangle \cong
(K_{n-1}/I_{n-1}) [y_n; \overline \tau_n],$ and thus, $K_n / \langle T
\rangle $ is a domain.  

If $y_n \in T$ but $x_n \notin T,$ then $K_n/\langle T \rangle \cong
(K_{n-1}/I_{n-1})[x_n; \overline \sigma_n]$ is a domain.  

If $x_n, y_n \in T,$ then $K_n/ \langle T \rangle \cong 
K_{n-1}/I_{n-1},$ and $K_n/ \langle T \rangle$ is a domain.  

If $\Omega_n \in T$ but $x_n, y_n \notin T,$ then $K_n / \langle T
\rangle \cong ((K_{n-1}/I_{n-1})[x_n; \overline \sigma_n][y_n; \overline
\tau_n, \overline \delta_n])/\langle \overline
\Omega_n \rangle.$  Since $T$ is admissible, $\Omega_{n-1} \notin T,$ so
 $\Omega_{n-1} \notin \la T_{n-1} \ra,$ and hence, $\overline
\Omega_{n-1}$ is a non-zero element of $K_{n-1}/I_{n-1}.$  Setting
$\Omega = q_n (q_n - p_n)^{-1} \overline \Omega_n,$ we have that $\langle
\Omega \rangle = \langle \overline \Omega_n \rangle.$  Then $\Omega =
x_ny_n + p_n (q_n - p_n)^{-1} \overline \Omega_{n-1},$ and $K_n / 
\langle T \rangle \cong ((K_{n-1}/I_{n-1})[x_n; \overline \sigma_n][y_n;
\overline \tau_n, \overline \delta_n])/\langle \Omega \rangle$ is a
domain by Lemma 2.6.  

We now have that $\langle T \rangle$ is completely prime in all cases.  
Since every admissible set $T$ consists of $H$-eigenvectors, it follows
that $\la T \ra$ is also $H$-prime. $\blacksquare$

\head
{3. $H$-prime ideals and $H$-simple localizations}
\endhead

Through the use of localizations, we will show that every $H$-prime
ideal of $K_n$ is generated by an admissible set.

\proclaim {Lemma 3.1}  Let $R$ be a $\bold k$-algebra and a domain, and
let $\sigma$ be a $\bold k$-automorphism of $R.$  Suppose that $H$ acts
on $R[x^{\pm1}; \sigma]$ so that $x$ is an $H$-eigenvector and $R$ is
both $H$-stable and $H$-simple, where $H$ acts on $R$ by restriction. 
If $H$ contains an automorphism $f$ such that $f \mid _R = \sigma,$ and
if $f(x) = \beta x$ for some $\beta \in \bold k^{\times}$, where $\beta$
is not a root of unity, then $R[x^{\pm1}; \sigma]$ is $H$-simple. 
\endproclaim

\noindent $\bold {Proof:}$ Let $I$ be a proper nonzero $H$-ideal of
$R[x^{\pm 1};\sigma],$ and note that, since $R$ is $H$-simple, $R \cap I
= 0.$  Then there exists $a \in I,$ with $a \neq 0,$ of shortest length
with respect to $x,$ say $a = a_{\ell}x^{\ell}+ \dots + a_m x^m$ for
some $\ell \leq m,$ where $a_i \in R$ for each $i$ and $a_m, a_{\ell}
\neq 0.$  Now, $x$ is a unit, so without loss of generality, $\ell = 0$
and $a = a_0 + a_1 + \dots + a_mx^m,$ where $m >0$ since $a \notin R.$   

Set $J = \{r \in R \mid r + r_1x + \dots + r_mx^m \in I $ for some 
$r_1, \dots r_m \in R \}$ and note that $J$ is an ideal of $R.$  Given
any $h \in H,$ let $\lambda_h$ be the $H$-eigenvalue of $x.$  Since $I$
is $H$-stable, we have that $h(r + r_1 x + \dots + r_mx^m) = h(r) +
h(r_1)\lambda_h x + \dots +h(r_m) \lambda_h ^m x^m \in I,$ so $h(J)
\subseteq J.$  Analogously, $h^{-1}(J) \subseteq J,$ so $h(h^{-1}(J)) = J
\subseteq h(J),$ yielding that $h(J) = J$ for each $h \in H.$  Hence,
$J$ is an $H$-stable ideal of $R.$

Now, $R$ is $H$-simple, so either $J = 0$ or $J = R;$ by our
choice of $a,$ we have that $J \neq 0,$ so $J = R.$  Then $1 \in J,$
and, without loss of generality, $a = 1 + a_1x + \dots + a_mx^m.$  Since
$I$ is $H$-stable, $f(a) = 1 + \sigma(a_1)\beta x + \dots + \sigma(a_m)
\beta^m x^m \in I,$ with $$f(a) - a = (\sigma(a_1)\beta - a_1)x + \dots
+ (\sigma(a_m) \beta^m - a_m)x^m \in I,$$ so that length $(f(a)-a ) \leq
m-1.$  Then $f(a) -a = 0$ and $a_i = \sigma(a_i)\beta^i$ for each $i.$ 
 
Consider $xax^{-1} = 1 + xa_1 + \dots + xa_mx^{m-1} = 1 + \sigma(a_1) x +
\dots + \sigma(a_m)x^m,$ and note that $xax^{-1} - a = (\sigma(a_1)-a_1)x
+ \dots + (\sigma(a_m)-a_m)x^m \in I.$  Then $\sigma (a_i) = a_i$ for
each $i,$ and in particular, $\sigma(a_m) = a_m = \sigma(a_m)\beta^m,$
so that $\beta^m = 1,$ a contradiction.  As a result, $R[x^{\pm1};
\sigma]$ contains no proper $H$-ideals.  
$\blacksquare$

\smallskip
  
\proclaim{Lemma 3.2}   Let $A = B[x; \sigma][y; \tau],$ where $B$ is a 
noetherian $\bold k$-algebra and both $\sigma$ and $\tau$ are $\bold
k$-automorphisms, such that $\tau(B) = B$ and $\tau(x) = \alpha x$ for
some $\alpha \in \bold k^{\times}.$    Suppose further that $B$ is a
domain, and that $H$ is a group of $\bold k$-automorphisms of $A$ such
that $B$ is $H$-stable and  $x,y$ are $H$-eigenvectors.   If there exist
$f,g \in H$ such that $f \mid_B = \sigma$ with $f(x) = \beta x$ and $g
\mid_{B[x; \sigma]} = \tau$ with $g(y) = \eta  y$ for some $\beta, \eta
\in \bold k^{\times},$ where $\beta, \eta$ are not roots of unity,  and
if $B$ is $H$-simple, then:

\rm{(a)  $A[x^{-1}][y^{-1}],$ $A/ \langle x,y \rangle,$ $(A/\langle x
\rangle)[y^{-1}],$ and $(A/ \langle y \rangle )[x^{-1}]$  are
$H$-simple.}

\rm{ (b) $A$ has precisely four $H$-prime ideals, namely:  $0,$
$\langle x \rangle,$ $\langle y \rangle,$ and $\langle x,y \rangle.$} 
\endproclaim

\smallskip

\noindent $\bold {Proof:}$  (a)  As in the proof of Lemma 2.6, $A[x^{-1}]
= B[x^{\pm1}; \sigma][y; \tau],$ and $A[x^{-1}][y^{-1}] = B[x^{\pm1};
\sigma][y^{\pm1}; \tau].$  By Lemma 3.1, $B[x^{\pm1}; \sigma]$ is
$H$-simple.  Noting that $B[x^{-1}; \sigma]$ is a domain as
a skew-Laurent ring over a domain, and that $g \mid_{B[x^{\pm1}; \sigma]}
= \tau,$ we apply Lemma 3.1 a second time to obtain that $B[x^{\pm1};
\sigma][y^{\pm1}; \tau] = A[x^{-1}][y^{-1}]$ is $H$-simple.

Since $A/ \langle x, y \rangle \cong_H B,$ it follows that $A/ \langle
x, y \rangle$ is $H$-simple.  Next, $(B[y; \tau])[y^{-1}] = B[y^{\pm1};
\tau]$ yields that $(A/ \langle x \rangle)[y^{-1}] \cong_H B[y^{\pm1};
\tau].$  By Lemma 2.6, $B[y^{\pm1}; \tau]$ is $H$-simple and hence,
$(A/\langle x \rangle)[y^{-1}]$ is $H$-simple.  Analogously,
$(A/ \langle y \rangle) [x^{-1}] \cong_H B[x^{\pm1}; \sigma]$ is
$H$-simple.  

(b)  Since $A = B[x; \sigma][y; \tau]$ is a domain, $0$ is an $H$-prime
ideal of $A.$  Observe that $A/ \la x \ra,$ $A/\la y \ra,$ and
$A/\la x,y\ra$ are domains since $A/\la x \ra \cong B[y; \tau],$
$A/ \la y \ra\cong B[x; \sigma],$ and $A/ \la x, y \ra \cong B.$ 
Further, $\la x \ra,$ $\la y \ra,$ and $\la x,y \ra$ are $H$-stable as
ideals generated by $H$-eigenvectors, so  $\la x \ra,$ $\la y \ra,$ and
$\la x,y \ra$ are $H$-prime.

Suppose that $P$ is an $H$-prime ideal of $A$ that does not contain
either $x$ or $y.$  Since $A$ is noetherian, $A[x^{-1}][y^{-1}]$ is
noetherian by \cite{\GW, Corollary 9.18}.  By Observation 1.9, $P$
extends to an $H$-stable prime ideal $Q = P^e = PA[x^{-1}][y^{-1}].$  By
(a), $A[x^{-1}][y^{-1}]$ is $H$-simple, so $Q = 0.$  Hence, $0 = Q^c =
\{q \in A \mid q1^{-1} \in Q \},$ and since $P \subseteq P^{ec} = Q^c
\subseteq 0,$ we have that $P = 0.$  Consequently, each nonzero
$H$-prime ideal of $A$ contains $x$ or $y.$ 

Now let $P$ be an arbitrary nonzero $H$-prime ideal and note that, for
any $i\in \Bbb{Z}^+,$ if $y^i \in P,$ then $y^{i-1}Ay = y^{i-1} \tau(A)y 
=y^iA \subseteq P.$  By Observation 1.8, either $y^{i-1} \in P$ or $y \in
P.$  Repeated applications of  Observation 1.8 in the first case yield
that $y \in P.$  Thus, if $x \in P$ and $y \notin P,$ then
$y^i \notin P$ for all $i \in \Bbb{Z}^+.$  Then $\overline {y^i} \notin P/ \la x \ra$ for each $i,$
and by Observation 1.9, $P/ \la x \ra$ corresponds to an $H$-prime ideal
$Q$ of $(A/\la x \ra)[y^{-1}].$  Since $(A/\la x \ra)[y^{-1}]$ is
$H$-simple, $Q = 0$ in $(A/\la x \ra)[y^{-1}],$ and $P = \la x \ra$
because $A$ is a domain.

Similarly, if $x^i \in P,$ then $x^{i-1}Ax = x^{i-1}\sigma(A)x = x^iA \subseteq P$ yields that $x \in P.$
As a result, whenever $y \in P$ with $x \notin P,$ we have that $P/\la y \ra$ corresponds to an 
$H$-prime ideal $I$ of $(A/\la y \ra)[x^{-1}].$  Since $(A/\la y \ra)[x^{-1}]$ is $H$-simple, 
$I = 0,$ so $P = \la y \ra.$

Finally, if $x \in P$ and $y \in P,$ then $J = P/ \la x,y \ra$ is
an $H$-prime ideal of $(A/\la x,y \ra)$.  Since $(A/\la x, y \ra)$ is
$H$-simple, $J=0,$ or $P = \la x, y \ra.$  We conclude that $0,$ $\la x
\ra,$ $\la y \ra,$ and $\la x, y \ra$ are the only $H$-prime ideals of
$A$. $\blacksquare$

\smallskip

\definition{Definition 3.3}  For $Q \in \spec K_n, $ let 
$$(Q:H)= \bigcap_{h \in H}h(Q).$$
\enddefinition

Note that $(Q:H)$ is an $H$-prime ideal.

\proclaim{Lemma 3.4} Let $A = B[x; \sigma][y; \tau, \delta],$
where $B$ is a noetherian $\bold k$-algebra and $\sigma,$
$\tau$ are both $\bold k$-automorphisms with $\tau(B) = B$
and $\tau(x) = \alpha x$ for some $\alpha \in \bold
k^{\times}.$ Assume that $\delta \mid_B = 0$ and that
$\delta(x) \in B$ is both nonzero and normal in $A.$ Suppose
further that $B$ is a domain, and that there exists an
element of the form $\Omega = xy + \lambda \delta(x),$ where
$\lambda \in \k ^{\times},$ such that $\Omega $ normalizes
$B[x; \sigma]$ and is normal in $A.$ Let $H$ be a group of
$\bold k$-automorphisms of $A$ such that $B$ is $H$-stable
and $x,$ $y,$ and $\Om$ are $H$-eigenvectors. Suppose that
there exist $f, g \in H$ such that $f \mid_B = \sigma$ with
$f(x) = \beta x$ and $f(y)= \hat \beta y$ such that $\alpha
\beta \hat \beta$ is not a root of unity and $g
\mid_{B[x^{\pm1};\sigma]} = \tau$ with $g(x^{-1}\Omega) =
\gamma (x^{- 1}\Omega)$ for some $\beta,$ $\hat \beta,$
$\gamma$ $\in \k^{\times},$ where $\beta,$ $\hat \beta,$ and
$\gamma$ are not roots of unity. If $B$ is $H$-simple, then:

{\rm (a)} $\delta(x)$ is invertible in A. 

{\rm (b)} no proper $H$-stable ideal of $A$ can contain a
power of $x.$ 

{\rm (c)} $A[x^{-1}][\Om^{-1}],$ $A[\Om^{-1}],$ and $(A/
\langle \Om\rangle)$ are $H$-simple. 

{\rm (d)} the only $H$-prime ideals of $A$ are $0$ and
$\langle \Om \rangle$.
\endproclaim 

\noindent $\bold {Proof:}$ (a) Since $\delta(x)
= yx - \alpha xy$ is an $H$-eigenvector, $\la \delta(x) \ra$
is an $H$-stable ideal of $A.$ Let $I$ be the ideal of $B$
generated by $\delta(x),$ and note that $\delta(x) \neq 0$
implies that $I \neq 0.$ Then $I$ is a nonzero $H$-stable
ideal of an $H$-simple ring, and hence, $I = B.$ In
particular, $1 \in I \subseteq \la \delta(x) \ra,$ so $\la
\delta(x) \ra = A.$ Since $\delta(x)$ is normal, $\delta(x)A
= A\delta(x) = \la \delta(x) \ra = A.$ Consequently,
$\delta(x)$ is invertible in $A.$

(b) Suppose that $P$ is a proper $H$-ideal of $A$ such that $x^j \in P$
for some $j > 0.$ Note that $x \notin P$ since $yx - \alpha xy
=\delta(x)$ and $\delta(x)$ is invertible by (a). Whenever $x^j\in P$
for some $j >1,$ we have that
$\delta(x^j) = yx^j -\tau(x^j)y = yx^j - \alpha^j x^jy \in
P.$ Further, $\sigma( \delta(x)) = f(\delta(x)) = f(yx -
\alpha xy) = \hat \beta \beta (yx-\alpha xy) = \hat \beta
\beta \delta(x).$ By \cite{\G, Lemma 1.1}, $$\align
\delta(x^j) &= \sum_{t = 0}^{j-1} [\tau(x)]^{t} \delta(x)
x^{j-1-t} = \sum_{t = 0}^{j-1} \alpha^{t} x^{t} \delta(x)
x^{j-1-t}\\ &= \sum_{t = 0}^{j-1} \alpha^{t}
\sigma^{t}(\delta(x)) x^{j-1} = \sum_{t=0}^{j-1} \alpha^t
(\hat \beta \beta)^t \delta(x) x^{j-1} = \Bigl(\sum_{t=0}^{j-1}
(\alpha \hat \beta \beta)^t \Bigr) \delta(x) x^{j-1}. \endalign $$
Recalling that $\alpha \hat \beta \beta$ is not a root of
unity, we have that $$\sum^{j-1}_{t=0} (\alpha \hat \beta
\beta)^t \neq 0, $$ and hence, $\delta(x) x^{j-1} \in P.$
Since $\delta(x)$ is invertible in $A,$ it follows that
$x^{j-1} \in P,$ and repeated applications of the above
argument yield that $x \in P,$ a contradiction. Thus, no
proper $H$-ideal of $A$ contains a power of $x.$ 

(c) By Lemma 2.6, $A[x^{-1}] = B[x^{\pm1};\sigma][x^{-1}\Omega;\tau],$
and hence,
$$A[x^{-1}][\Omega^{-1}]
=B[x^{\pm1};\sigma][x^{-1}\Omega;\tau][\Omega^{-1}].$$
Let $C=B[x^{\pm1};\sigma][x^{-1}\Omega; \tau][\Omega^{-1}],$ and
note that $C = B[x^{\pm1};\sigma][(x^{-1}\Omega)^{\pm1};
\tau].$ Applying Lemma 3.1 twice yields that both $B[x^{\pm1};
\sigma]$ and $C$ are $H$-simple, so $A[x^{-1}][\Om^{-1}]$ is
$H$-simple. 

Let $P$ be an $H$-prime ideal of $A[\Om^{-1}].$  By
Observation 1.9, $P$ is induced from  an
$H$-prime ideal $\hat P$ of $A$ disjoint from $\{\Om^i \mid i
\in \Bbb{Z^+} \}.$  By (b), $\hat P$ is also disjoint from
$\{x^j \mid j \in \Bbb{Z^+} \}.$  

Suppose that there exist some $i, j\in \Bbb{Z}$$_{\geq 0}$
such that $\Om^ix^j \in \hat P.$  Since $\Om$ is normal in $A,$ 
we have that $\Om^i A x^j = A \Om^i x^j \subseteq \hat P$ where 
$\Om^i, x^j$ are $H$-eigenvectors. By Observation 1.8, it follows that
either $\Om^i \in \hat P$ or $x^j \in \hat P,$ a contradiction.  Thus,
$\hat P$ is disjoint from the multiplicative set generated by $x$ and
$\Om.$  

By Observation 1.9, $\hat P$ corresponds to an $H$-prime
ideal $\tilde P$ of $A[x^{-1}][\Om^{-1}].$  Since this
algebra is $H$-simple, $\tilde P = 0,$ and hence, $\hat P
=0,$ so $P = 0.$  Thus, $A[\Om^{-1}]$ contains no nonzero
$H$-prime ideals.

If $I$ is a proper $H$-ideal of $A[\Om^{-1}],$ then $I$ is contained 
in some prime ideal $J$ of $A[\Om^{-1}].$  Further, $Q = (J:H)$ is an 
$H$-prime ideal with $I \subset J \subset Q.$  Since $A[\Om^{-1}]$
contains no nonzero $H$-prime ideals, it follows that $Q = I = 0.$ 
Thus, $A[\Om^{-1}]$ is $H$-simple.   

Next, note that $\la \Om \ra = A \Om
= \Om A$ since $\Om$ is normal in $A,$ and by \cite{\GW,
Theorem 9.20(a)}, $A[x^{-1}](A \Om) = (A \Om)A[x^{-1}]$ is a
two-sided ideal of $A[x^{-1}].$ Then $A[x^{-1}] \Om
=A[x^{-1}] (A \Om)$ and $(A \Om) A[x^{-1}] = (\Om A)
A[x^{-1}] = \Om A[x^{-1}]$ yield that $A[x^{-1}] \Om = \Om
A[x^{-1}].$ Hence, the ideal generated by $\Omega$ in $A[x^{-1}]$ is
$\Om A[x^{-1}],$ and
$$A[x^{-1}]/(\Omega A[x^{-1}]) = (B[x^{\pm1}; \sigma][x^{-1}\Omega;
\tau]) /(\Omega B [x^{\pm1}; \sigma][x^{-1}\Omega; \tau]) \cong_H
B[x^{\pm1};\sigma]$$ via the map sending $ \sum (b_ix^i +
\Omega A[x^{-1}])$ to $\sum b_ix^i.$ Thus,
$A[x^{-1}]/(\Omega A[x^{-1}])$ is $H$-simple. Now, $(A /
\langle \Omega \rangle)[x^{-1}] \cong_H A[x^{-1}]/(\Omega
A[x^{-1}]),$ and hence, $(A/ \langle \Omega
\rangle)[x^{-1}]$ is $H$-simple. Let $w = (\delta(x))^{-1}$
in $A,$ and consider the product $\overline x \overline y$
in $A/ \la \Om \ra.$ We have that $\overline x \overline y =
-\lambda \overline{\delta(x)},$ so $-\lambda^{-1} \overline x
\overline y \overline w = \overline{\delta(x) w} = \overline
1.$ Thus, $\overline x$ has a right inverse in $A/ \la \Om
\ra$ and $\overline x$ is invertible since $A$ is a domain.
As a result,$ A/ \la \Om \ra = (A/ \langle \Omega
\rangle)[x^{-1}]$ is $H$-simple. 

(d) Note that $0$ is an $H$-prime ideal of $A$ since $A$ is a domain.
Further, $\la \Om \ra$ is $H$-stable because $\Om$ is an $H$-eigenvector,
and $A/\la \Om \ra$ is a domain by Lemma 2.6, so $\la \Om \ra$
is a (completely) prime $H$-ideal of $A.$ Now let $P$ be any
$H$-prime ideal of $A.$ If $P$ is disjoint from $S= \{c
\Omega^i \mid c \in \k^{\times} $ and $i \in \Bbb{Z}$$_{\geq 0}\},$ then
$P$ extends to an $H$-prime ideal $\widehat P$ of $A[\Omega^{-1}]$ by
Observation 1.9. By (c), the localization $A[\Om^{-1}]$ is
$H$-simple, so $\widehat P =0,$ and $P = 0.$ Now assume that
there exists some $c\Om^i \in P \cap S.$  Then $(\Om A)^i = \Om^iA
\subseteq P,$  and successive applications of Observation 1.8 yield that
$\Om \in P.$ Hence, $P/ \la \Om \ra$ is an $H$-ideal of $A/ \la \Om \ra.$
By (c), $A/ \la \Om \ra$ is $H$-simple, so $P = \la \Om \ra.$ Thus, $0$
and $\la \Om \ra$ are indeed the only $H$-prime ideals of $A.$
$\blacksquare$

\smallskip 
Recall the set $N_T$ defined in Definition 2.3.

\proclaim {Lemma 3.5}  Let $T$ be an admissible set.  Then each element of
$N_T$ represents a nonzero,  normal coset in $K_n/ \la T \ra.$
\endproclaim

\smallskip

\noindent {\bf Proof:}  By Theorem 2.8, $\la T \ra \cap \P_n = T,$ so each
element of $N_T$ represents a nonzero coset in $K_n/ \la T \ra.$   Since
$\Om_i$ is normal for each $i,$ we have that $\overline \Omega_i$ is
normal and nonzero whenever $\Om_i \in N_T.$  Similarly, $\overline x_1$
and $\overline y_1$ are nonzero, normal elements of $K_n / \la T \ra$
whenever $x_1 \in N_T$ and $y_1 \in N_T,$  respectively.

Let $i > 1$ and suppose that $x_i \in N_T.$  As noted above, $\overline
x_i$ is nonzero.  Further, $x_i \notin T$ with either $\Omega_{i-1} \in
T$ or $\Omega_i \in T,$ so $\overline x_i$ is normal in 
$K_n / \la T \ra$ by Lemma 2.1.  Analogously, $\overline y_i$ is normal
and nonzero whenever $y_i \in N_T.$ 
$\blacksquare$

\smallskip

\definition {Definition 3.6}  Given an admissible set $T,$ let $E_T$ be
the multiplicative set generated by $N_T \cup \k^{\times}.$  Then $E_T$
is
$H$-stable and $\overline E_T$ forms a denominator set of $K_n / \la T
\ra.$ 
\enddefinition                

\smallskip

\proclaim{Theorem 3.7} Given an admissible set $T,$ let $K_n^T =
(K_n/\la T \ra)[E_T^{-1}].$  Then $K_n^T$ is $H$-simple.
\endproclaim

\noindent $\bold {Proof:}$  We will proceed by induction on $n,$ and 
for $n = 1,$ we will use Lemma 3.2.  By Proposition 2.5, 
$K_1 = \k[x_1][y_1; \tau_1],$ so in the format of Lemma 3.2, $\sigma$
is the identity map.  Let $\beta \in \k^{\times}$ such that $\beta$ is
not a root of unity and consider $f = (\beta, 1) \in H.$  Then $f$ acts
as the identity on $\k$ with $f(x_1)= \beta x_1.$  Further, let $g =
(q_1^{-1}, \beta).$  Then $g$ acts as $\tau_1$ on $\k[x_1]$ with 
$g(y_1) = \beta y_1.$   

Now, the four possible cases for $T$ are: $ \emptyset, $ 
$\{x_1, \Om_1 \},$ $\{y_1, \Om_1 \},$ and  $\P_1.$  As a result, 
$K_1^T= K_1[x_1^{-1}][y_1^{-1}],$ $(K_1/\la x_1\ra)[y_1^{-1}],$ 
$(K_1/\la y_1\ra)[x_1^{-1}],$ or $K_1/\la x_1, y_1\ra.$  Applying 
Lemma 3.2, yields that each of the possible cases for $K_1^T$ is
$H$-simple.

Suppose now that $n > 1$ and $K_{n-1}^S$ is $H$-simple for any
admissible set $S \subseteq \P_{n-1}.$  Given an  admissible set $T$ of
$K_n,$ set $T_{n-1}=T \cap \P_{n-1},$ and let $I_{n-1}$ be the ideal of
$K_{n-1}$ generated by $T_{n-1}.$  Then, as in the proof of Theorem 2.8,
$$K_n/ \la T_{n-1} \ra \cong_H (K_{n-1}/I_{n-1}) [\overline x_n;\overline {\sigma}_n] [\overline y_n; \overline                
{\tau}_n; \overline {\delta}_n], $$ where $\overline {\delta}_n = 0$ if
$\Om_{n-1} \in T_{n-1}.$  Note that $\tau_n(E_{T_{n-1}}) =
E_{T_{n-1}}$ and $\sigma_n( E_{T_{n-1}}) = E_{T_{n-1}}.$  Then applying 
\cite{\G, Lemma 1.4} twice yields that $$(K_n/ \la T_{n-1} \ra)
[E_{T_{n-1}}^{-1}] \cong_H (K_{n-1}/I_{n-1})[E_{T_{n-1}}^{-1}] 
[\overline x_n;\overline {\sigma}_n] [\overline y_n; \overline {\tau}_n;
\overline {\delta}_n]. $$ Let $R= (K_{n-1}/I_{n-1})[E_{T_{n-1}}^{-1}].$  

Setting $S = T \setminus T_{n-1},$ note that $\la T \ra = \la T_{n-1}
\ra + \la S \ra.$  Then $K_n/\la T \ra \cong_H (K_n/\la T_{n-1}
\ra)/\la S \ra$, and $(K_n/ \la T \ra )[E_{T_{n-1}}^{-1}] \cong_H
(R[\overline x_n; \overline \sigma_n][\overline y_n; \overline \tau_n,
\overline \delta_n])/ \la S \ra.$  Define $E$ to be the multiplicative
set generated by $N_T \setminus(N_{T_{n-1}} \cap P_{n-1}).$  Then $E$
forms a denominator set for $K_n/ \la T \ra$ with
$(\K/t)[E_T^{-1}] = (\K/t) [E_{T_{n-1}}^{-1}][E^{-1}] \cong_H
((R[\overline x_n; \overline \sigma_n][\overline y_n; \overline
\tau_n, \overline \delta_n])/ \la S \ra)[E^{-1}].$   

In order to apply Lemma 3.2 and Lemma 3.4 below, we will first define the
necessary elements of $H.$  Let $\beta \in \k^{\times}$ such that $\beta$
and $\hat \beta = \beta^{-1}p_n$ are not roots of unity.  Next, set
$$\xalignat 1 \hat f & = (q_1^{-1}p_n \gamma_{n,1}, q_1\gamma_{1,n},
\dots, q_{n-1}^{-1}p_n \gamma_{n, n-1}, q_{n-1}\gamma_{n-1, n}, \beta,
\beta^{-1}p_n)
\text { and } \\ 
\hat g &= (p_n^{-1}\gamma_{1,n}, \gamma_{n,1}, \dots, p_n^{-1}
\gamma_{n-1,n}, \gamma_{n,n-1}, q_n^{-1}, p_n^{-1}q_n).
\endxalignat$$ 
Then $\hat f, \hat g \in H$ with $\hat f \mid_{K_{n-1}} = \sigma_n,$
$\hat f(x_n) = \beta x_n,$ and $\hat f(y_n) = \hat \beta y_n.$ 
Also, $\hat g \mid_{K_{n-1}[x_n; \sigma_n]} = \tau_n,$ $\hat g(y_n) 
= \eta y_n$ where $\eta = p_n^{-1}q_n,$ and $\hat g(x_n^{-1} \Om_n)
= \eta x_n^{-1}\Om_n.$  As defined, $\beta, \eta, \hat \beta,$ and
$q_n^{-1}\beta \hat \beta = q_n^{-1}p_n$ are not roots of unity.      

The five possible cases for $S$ are:  $\emptyset,$ $\{x_n, \Om_n \},$
$\{y_n, \Om_n\},$ $\{\Om_n\},$ and $\{x_n, y_n, \Om_n\}.$
If $S =\emptyset,$ then $\la S \ra = 0,$ and if $\Om_{n-1} \in T_{n-1},$
then $E$ is generated by $x_n$ and $ y_n,$ so that $K_n^T \cong_H
((R[\overline x_n;\overline \sigma_n][\overline y_n; \overline \tau_n,
\overline\delta_n])/ \la S \ra)[E^{-1}] = (R[\overline x_n; \overline
\sigma_n][\overline y_n; \overline \tau_n, \overline
\delta_n])[x_n^{-1}][y_n^{-1}].$  Applying Lemma 3.2 yields that
$K_n^T$ is $H$-simple.  If $\Om_{n-1} \notin T_{n-1},$ then 
$E$ is generated by $\Om_n$ and $K_n^T \cong_H (R[\overline x_n;
\overline \sigma_n][\overline y_n; \overline \tau_n, \overline
\delta_n])[\Om_n^{-1}]$ is $H$-simple by Lemma 3.4.

If $S = \{x_n, \Om_n \},$ then $\la S \ra = \la x_n \ra$ and $E$ is
generated by $y_n.$  Further, $x_n \in T$ implies that $\Om_{n-1} \in
T_{n-1},$ so $\overline {\delta}_n = 0$ and $(R[\overline x_n;
\overline \sigma_n][\overline y_n; \overline \tau_n,
\overline \delta_n]) = (K_{n-1}/I_{n-1})[E_{T_{n-1}}^{-1}][\overline
x_n; \overline \sigma_n][\overline y_n; \overline \tau_n].$ By Lemma 3.2,
$K_n^T \cong_H ((R[\overline x_n; \overline \sigma_n][\overline y_n;
\overline \tau_n])/\la x_n \ra)[y_n^{-1}]$ is $H$-simple. 

If $S = \{y_n, \Om_n \},$ then $\la S \ra = \la y_n \ra$ and $E$ is
generated by $x_n.$  Moreover, $\overline \delta_n = 0,$ so $K_n^T
\cong_H ((R[\overline x_n; \overline \sigma_n][\overline y_n; \overline
\tau_n])/ \la y_n \ra)[x_n^{-1}]$ is $H$-simple by
Lemma 3.2.

If $S = \{ \Om_n \},$ then $\Om_{n-1} \notin T_{n-1}$ with
$E=\k^{\times},$ and $K_n^T \cong_H (R[\overline x_n; \overline
\sigma_n] [\overline y_n; \overline \tau_n, \overline \delta_n])/
\la \Om_n \ra$ is $H$-simple by Lemma 3.4.

Lastly, if $S = \{x_n, y_n, \Om_n \},$ then $\Om_{n-1} \in T$ and 
$E = \k^{\times}.$  By Lemma 3.2, $K_n^T \cong_H (R[\overline x_n;
\overline \sigma_n][\overline y_n; \overline \tau_n])/ \la x_n, y_n
\ra$ is $H$-simple.  We conclude that $K_n^T$ is $H$-simple for any
admissible set $T.$ 
$\blacksquare$    

\smallskip

\proclaim{Theorem 3.8}  Every $H$-prime ideal of $K_n$ is generated by an
admissible set.
\endproclaim

\noindent $\bold {Proof:}$  Let $P$ be an $H$-prime ideal of $K_n,$ and
let $T = P \cap \P_n,$ noting that $T$ is admissible by Lemma 2.4.  By
definition, $N_T \cap T = \emptyset,$ so $N_T \cap P = \emptyset.$ 
Now, $P/ \la T \ra$ is an $H$-prime ideal of $K_n/ \la T \ra$ with
$\overline N_T  \cap P/ \la T \ra = \emptyset,$ where each element of 
$\overline N_T$ is normal in $K_n /\la T \ra.$  Recalling that 
$\overline E_T$ is the multiplicative set generated by $N_T,$ we have
that $\overline E_T \cap P/ \la T \ra = \emptyset$ by Observation 1.8. 
Since $K_n / \la T \ra$ is a domain by Theorem 2.8, Observation 1.9 yields
that $(P/ \la T \ra)[E_T^{-1}]$ is an $H$-prime ideal  of $K_n^T.$  By
Theorem 3.7, $K_n^T$ is $H$-simple, yielding that $(P/ \la T \ra)[E_T^{-1}]
= 0$ in $K_n^T.$   Then $P/ \la T \ra = 0$ by Observation 1.9, so $P =
\la T \ra.$ $\blacksquare$        

\proclaim{Corollary 3.9}  Every $H$-prime ideal of $K_n$ is completely
prime. 
\endproclaim

\noindent $\bold {Proof:}$  By Theorem 2.8, every ideal generated by an
admissible set is a completely prime ideal.
$\blacksquare$

\head
{4.  Primitive Ideals and Catenarity}
\endhead

We will show how the $H$-prime ideals of $K_n$ are related to the
algebra's primitive ideals.  As preliminary steps, we will see that
$K_n$ satisfies the Nullstellensatz over $\bold k$ and is normally 
separated.  Lastly, we will show that $K_n$ is catenary.

\definition{Definition 4.1}  Let $\Hspec K_n$ denote the set of $H$-prime
ideals of $K_n.$ 
\enddefinition

\definition {Definition 4.2}  Let $G$ be an affine algebraic group over
$\k$ and let $G$ act on a $\k$-algebra $A$  by $\k$-algebra
automorphisms.  Then $G$ acts {\it rationally \/} on $A$ if $A$ is a
directed union of finite  dimensional $G$-invariant $\k$-subspaces $V_i$
with the property that each of the restriction maps 
$G \rightarrow \Aut A$ $\rightarrow GL(V_i)$ is a morphism of
algebraic varieties.
\enddefinition

This definition simplifies considerably when $G$ is an algebraic
torus.  See \cite{\BG, II.2.6}.

\proclaim {Theorem 4.3}  If $\k$ is an infinite field, then $H$ acts
rationally on $K_n.$  
\endproclaim

\noindent $\bold {Proof:}$  By Proposition 2.5, $K_n$ has a $\k$-basis  
$$ \Cal {A} = \{x_1^{r_1}y_1^{r_2}...x_n^{r_{2n-1}}y_n^{r_{2n}} 
\mid r_i \in \Bbb {Z}^+ \}.$$  Since each $x_i$ and $y_i$ is an
$H$-eigenvector, $\Cal {A}$ consists of $H$-eigenvectors, so $K_n$ is a
direct sum of $H$-eigenspaces.  Fix $z = x_1^{r_1}y_1^{r_2} \dots
x_n^{r_{2n-1}}y_n^{r_{2n}}.$ For each $h \in H,$  we have that $h(z) =
(h_1^{r_1}h_2^{r_2} \dots h_{2n}^{r_{2n}})z$ where $h_1,\d, h_{2n} \in
\k^{\times}$ are the components of $H.$  Then the $H$-eigenvalue for $z$
is the character $f:H \longrightarrow \k^{\times}$ defined by
$f(h)=h_1^{r_1}h_2^{r_2} \dots h_{2n}^{r_{2n}}.$  

Letting $p_1, \d, p_{2n}:H \longrightarrow \k^{\times}$ be the
projection maps, we have that $p_1, \d, p_{2n}$ are rational characters
(see \cite{\BG, II.2.5}).  Then $f = p_1^{r_1}p_2^{r_2} \dots
p_{2n}^{r_n}$ is a rational character. By \cite{\BG, Theorem II.2.17},
$H$ acts rationally on $K_n.$ $\blacksquare$  

\definition {Definition 4.4}  Let $J$ be an $H$-prime ideal.  Then 
$$\spec_J K_n = \{P \in \spec K_n \mid (P:H) = J\}.$$  Note that all the
prime ideals in $\spec_J K_n$ contain $J.$ The set $\spec_J K_n$ is the
{\it H-stratum \/} of $\spec K_n$ correponding to $J.$   These
$H$-strata partition $\spec K_n:$ $$\spec K_n = \bigcup_{J \in \Hspec
K_n} \spec_J K_n.$$
    
\enddefinition

\proclaim{Theorem 4.5}  Let $\k$ be an infinite field and let $J \in
\Hspec K_n.$  Further, define $E_J$ to be the set  of all regular
$H$-eigenvectors in $K_n/J.$  Then:

{\rm (a)}  $E_J$ is a denominator set and the localization $K_n^J=
(K_n/J)[E_J^{-1}]$ is $H$-simple.  

{\rm (b)}  $\spec_J K_n$ is homeomorphic to $\spec K_n^J$ via
localization and contraction.

{\rm (c)}  $\spec K_n^J$ is homeomorphic to $\spec Z(K_n^J)$ via
contraction and extension.

{\rm (d)}  $Z(K_n^J),$ the center of $K_n^J,$ is a Laurent polynomial
ring, in at most $n+1$ indeterminates, over the fixed field $Z(K_n^J)^H
= Z($Fract $K_n/J)^H.$  The indeterminates can be chosen to be 
$H$-eigenvectors with linearly independent $H$-eigenvalues.  
\endproclaim

\smallskip

\noindent $\bold {Proof:}$  Note that $H \cong (\k^{\times})^{n+1}$ via
the mapping $$(h_1,h_2, \d, h_{2n-1}, h_{2n}) \longmapsto (h_1, h_2,
h_3, h_5, \d, h_{2n-3}, h_{2n-1})$$ because $h_{2i}=h_{2i-1}^{-1}h_1h_2$
for all $i=2, \d, n.$  Since $H$ acts rationally on $K_n$ by Theorem 4.3, 
the result follows from \cite{\BG, Theorem II.2.13}.
$\blacksquare$

\definition{Observation 4.6}  If $J = \la T \ra,$ then $K_n/J$ is a
domain by Theorem 2.8, so $\overline E_T \subseteq E_J$ and hence, $K_n^T$
is a subalgebra of $K_n^J.$
\enddefinition

\proclaim{Lemma 4.7}  If $T \subset S$ are admissible sets, then
$\overline S \cap \overline E_T \ne \emptyset$ in $\K/t.$  
\endproclaim

\noindent $\bold {Proof:}$  For each $j \in \{1, \d, n\},$ if $\Om_j
\in S $ with $\Om_j \notin T,$  then $\Om_j \in N_T \subset E_T,$ so 
$\overline \Om_j \in \overline E_T \cap \overline S.$ If every
$\Om_j \in S$ is also contained in $T,$ then either $x_i \in S \setminus
T$ or $y_i \in S \setminus T$ for some $i \in \{1, \d, n\}.$  
In the first case, $ \Om_i \in S$ since $S$ is admissible, so $x_i \in
N_T,$ and $\overline x_i \in \overline E_T \cap \overline S.$  In the
second case, $\overline y_i \in \overline E_T \cap \overline S$ by a 
similar argument.  $\blacksquare$     

\proclaim{Lemma 4.8}  If $\k$ is an infinite field and $J \in \Hspec K_n,$ let $T$ be the admissible set that generates $J.$  
Then there exists a one-to-one correspondence between $\spec K_n^T$ and $\spec K_n ^J$ given by extension and contraction.
\endproclaim

\noindent $\bold {Proof:}$  By \cite{\GW, Theorem 9.22}, there exists a bijection via extension and contraction
between $\spec K_n^T$ and the prime ideals of $K_n/ \la T \ra= K_n/J$ that are disjoint from $\overline E_T.$  
Since the ideals of $K_n/J$ are precisely the cosets of the ideals of $K_n$ that contain $J,$ we consider the 
set $$W = \{P \in \spec K_n \mid J \subseteq P \text{ and } \overline P \cap \overline E_T = \emptyset 
\text{ in } K_n / J  \}.$$  We will show that $W = \spec_J K_n,$ and the result will then follow from Theorem 4.5.  

Let $P \in \spec_J K_n$ and recall that $J \subseteq P.$  Then $P$ extends to a prime ideal 
of $K_n^J = (K_n/J)[E_J^{-1}]$ by Theorem 4.5, so $\overline P \cap E_J = \emptyset.$  By Observation 4.6,
we have that $\overline P \cap \overline E_T = \emptyset,$ yielding that $P \in W.$  Thus, $\spec_J K_n \subseteq W.$

Next, suppose that $Q \in W$ and consider $(Q:H).$  By Theorem 3.8, there exists an admissible set $S$ such that 
$\la S \ra = (Q:H).$  Note that $S = (Q:H) \cap P$ by Theorem 2.8, so $T \subseteq S$ with 
$\overline E_T \cap \overline S \subseteq \overline E_T \cap \overline Q = \emptyset.$  
Applying Lemma 4.7 yields that $S = T,$ or $\la T \ra = (P:H).$  Hence, $P \in \spec_J K_n,$ 
and $W \subseteq \spec_J K_n.$  Consequently, $W= \spec_J K_n$ as
desired. 
$\blacksquare$

\proclaim{Proposition 4.9}  Let $\k$ be an infinite field with $J \in 
\Hspec K_n$ generated by the admissible set $T.$  Then there exists a
one-to-one correspondence $\theta: \spec_J K_n \longrightarrow 
\spec Z(K_n^T)$ with $\theta(P) = (P/J)K_n^J \cap Z(K_n^T).$  Further,
letting $\psi_J: K_n \longrightarrow K_n/J \longrightarrow K_n^J$ be the
localization map, we have that $\theta^{-1}(Q)
= \psi_J^{-1}(QK_n^J)$ for each $Q \in \spec Z(K_n^T).$  Moreover, 
both $\theta$ and $\theta^{-1}$ preserve inclusion.
\endproclaim

\noindent {\bf Proof:}  Let $P \in \spec_J K_n.$  By Theorem 4.5, $P$
extends uniquely to $(P/J)K_n^J,$ a prime ideal of $K_n^J.$  Then by
Lemma 4.8, $(P/J)K_n^J$ contracts to a unique prime ideal, $(P/J)K_n^J
\cap K_n^T,$ of $K_n^T.$  Since $K_n^T$ is $H$-simple by Theorem 3.7,
contraction and extension provide bijections between $\spec K_n^T$
and $\spec Z(K_n^T)$ by \cite {\BG, Corollary II.3.9}.  As a result, we
then have  that $\Bigl( (P/J)K_n^J \cap K_n^T \Bigr) \cap Z(K_n^T) =
(P/J)K_n^J \cap Z(K_n^T) \in \spec Z(K_n^T).$  Consequently, the
map $\theta: \spec_J K_n \longrightarrow \spec Z(K_n^T),$ defined by
$\theta(P) = (P/J)K_n^J \cap Z(K_n^T),$ is indeed a bijection as the
composition of one-to-one correspondences.

To compute $\theta^{-1},$ let $Q \in \spec Z(K_n^T).$  Then $Q$ extends
to $QK_n^T \in \spec(K_n^T)$ by \cite{\BG, Corollary II.3.9}, and 
$QK_n^T$ extends to $QK_n^J \in \spec (K_n^J)$ by Lemma 4.8.  Then by
Theorem 4.5, $\psi_J^{-1}(QK_n^J) \in \spec_J K_n,$ so
$\theta^{-1}(Q) =\psi_J^{-1}(QK_n^J).$  Note that both $\theta$ and
$\theta^{-1}$ preserve inclusion by construction.
$\blacksquare$

\definition{Definition 4.10}  A ring $R$ is a {\it Jacobson ring \/} if
each prime ideal $P$ satisfies 
$J(R/P)=0.$  
\enddefinition

\definition{Definition 4.11} A noetherian $\k$-algebra $A$ {\it satisfies
the Nullstellensatz over $\k$ \/}  if:

(i) $A$ is a Jacobson ring.

(ii)  the endomorphism ring of every irreducible $A$-module is algebraic over $\k.$
\enddefinition

\proclaim{Theorem 4.12}  The algebra $K_n$ satisfies the Nullstellensatz over $\k.$
\endproclaim

\noindent $\bold {Proof:}$  Note that 
$$ \k = K_0 \subset K_0[x_1] \subset K_1 \subset K_1[x_2; \sigma_2] \subset K_2 \subset \dots \subset K_{n-1} 
\subset K_{n-1}[x_n; \sigma_n] \subset K_n$$
is a sequence of subalgebras of $K_n.$  For each $i \geq 0,$ the subalgebra $K_i[x_i; \sigma_i]$ is generated by
$K_i$ together with $x_i$ so that $(K_i)x_i = x_i(K_i).$  Further, for $i > 1,$ the subalgebra $K_i$ is generated 
by $K_{i-1}[x_i; \sigma_i]$ together with $y_i,$ satisfying 
$$(K_{i-1}[x_i; \sigma_i])y_i + K_{i-1}[x_i; \sigma_i] =y_i(K_{i-1}[x_i; \sigma_i]) + K_{i-1}[x_i; \sigma_i].$$  
By \cite{\MR,Theorem 9.4.21}, $K_n$ satisfies the Nullstellensatz over $\k$ 
(cf. \cite{\BG, Theorem II.7.17}).  $\blacksquare$

 \definition{Definitions 4.13}  Let $X$ be a topological space.  Then a
subset $C$ of $X$ is {\it locally closed \/} if there exists an open set
$U$ such that $C \subseteq U$ and $C$ is closed in $U.$   A point $x \in
X$ is a {\it locally closed point \/} if the singleton $\{ x \}$ is
locally closed.
\enddefinition

Given a ring $R,$ we say that $P \in \spec R$ is locally closed if $P$
is a locally closed point of  spec $R$ in the Zariski topology.  By
\cite{\BG, Lemma II.7.17}, this is equivalent to the condition  that
$\bigcap \{Q \in $ spec $R \mid P \subset Q \}$ is an ideal properly
containing $P.$  

\definition {Definition 4.14}  A prime ideal $P$ of a noetherian
$\k$-algebra $A$ is said to be {\it rational \/} if  the field $Z($
Fract $A/P)$ is algebraic over $\k.$  
\enddefinition

\definition {Definition 4.15}  An algebra $A$ satisifies the
{\it Dixmier-Moeglin \/} equivalence if the sets of primitive, rational,
and locally closed primes coincide.
\enddefinition

\proclaim {Theorem 4.16}  Let $\k$ be an infinite field. Then $K_n$
satisfies the Dixmier-Moeglin equivalence and  the primitive ideals of
$K_n$ are precisely the primes maximal in their $H$-strata:

$$\xalignat1 \prim K_n &= \{\text{locally closed prime ideals}\}\\
           &= \{\text{rational prime ideals}\}\\
           &= \cup_{J \in \Hspec K_n} \{\text{maximal elements of }
\spec_J K_n \}. \endxalignat$$
\endproclaim

\noindent $\bold {Proof:}$  By Theorem 4.3, the group $H$ acts
rationally on $K_n.$  Since each $H$-prime ideal is generated by an
admissible set, $\Hspec K_n$ is finite.  Further, $K_n$ satisfies  the
Nullstellensatz over $\k$ by Theorem 4.12 and \cite{\BG, Theorem
II.8.4} (a specialization of \cite{\GLet, Theorem 2.12}) yields the
desired results. $\blacksquare$

\proclaim {Corollary 4.17}  The primitive ideals of $K_n$ correspond
to maximal ideals of the various $Z(K_n^T)$ under the one-to-one
correspondences described in Proposition 4.9. 
\endproclaim

\definition {Definitions 4.18}  A chain of prime ideals 
$$P_0 \subset P_1 \subset P_2 \subset ... \subset P_{\ell}$$ 
 of a ring $R$ has {\it length \/} $\ell.$  The chain is {\it saturated
\/} if, for $1 \leq i \leq \ell,$ there is no prime ideal $P$ such that
$P_{i-1} \subset P \subset P_i.$  

If $Q$ is a prime ideal, the supremum of the lengths of all of the
chains of primes  contained in $Q$ is the {\it height \/} of $Q,$
denoted $ht (Q).$

The ring $R$ is {\it catenary \/} if for every pair of prime ideals $P$
and $Q$ of $R$ such that $P \subset Q,$ all saturated chains of primes
from $P$ to $Q$ have the same length.

Lastly, $\spec R$ is {\it normally separated \/} if, for every pair 
$P \subset Q$ of distinct primes of $R,$ there exists a nonzero element
of $Q \setminus P$ which is normal in $R/P.$

\enddefinition

\proclaim {Theorem 4.19}  The algebra $K_n$ is normally separated. 
\endproclaim

\noindent $\bold {Proof:}$  By \cite{\BG, Theorem II.9.15}, it suffices
to prove normal $H$-separation, meaning that for every pair of distinct
$H$-primes $I \subset J,$ the factor $J/I$ contains a nonzero normal 
$H$-eigenvector.  If $I \subset J$ are $H$-prime ideals of $K_n,$ then
applying Theorem 3.8 yields that $I = \la T_I \ra$ and $J = \la T_J \ra$
for some admissible sets $T_I \subset T_J.$  By Lemma 4.7,  there exists
some $z \in T_J$ such that $\overline z$ is normal and nonzero in $K_n /
\la T_I \ra  = K_n / I.$  Since each element of $\P_n$ is an
$H$-eigenvector, the result follows. $\blacksquare$ 

The {\it Gelfand-Kirillov dimension\/} of a $\k$-algebra $A$ is denoted
by $\GKdim(A).$  Further details may be found in \cite {\MR, Chapter 8}. 
The Auslander-regular and Cohen-Macaulay conditions are defined in
\cite{\BG, Appendix 1.5}, for instance. 

\proclaim {Theorem 4.20}  The algebra $K_n$ is catenary, and if $P$ and
$Q$ are prime ideals of $K_n$  such that $P \subset Q,$ then 
$$\height (Q/P) = \GKdim (K_n/P) - \GKdim (K_n/Q).$$
In particular, 
$$\height(Q) + \GKdim (K_n/Q) = 2n$$
for every $Q \in spec K_n.$
\endproclaim

\noindent $\bold {Proof:}$  By \cite{\BG, Lemma II.9.7}, $\GKdim(K_n) =
2n.$  Further, $K_n$ is Auslander-regular and Cohen-Macaulay by
\cite{\BG, Lemma II.9.10}.  Hence, we may apply \cite{\GLen, Theorem
1.6}. $\blacksquare$
\smallskip

We refer to the second formula of Theorem 4.20 by saying that {\it
Tauvel's height formula \/} holds in $K_n.$ 
     
\head
{5.  Some Examples}
\endhead

In this section, we will show how to compute the primitive ideals of
$K_2$ from the admissible sets.   Throughout, $\langle \dots \rangle$
will either denote an ideal of $K_n$ or a free abelian group  generated
by elements of $\bold k^{\times}.$
 
By Corollary 4.17, the primitive ideals of $K_n$ are those that
correspond to maximal ideals of the various $Z(K_n^T)$ under the
one-to-one correspondences of Proposition 4.9.  Thus, to find the
primitive ideals, we first determine $N_T$ and then localize $K_n / \la
T \ra$ with respect to the multiplicative set generated by the elements
of $N_T.$  Next, we find the generators of $Z(K_n^T);$ by \cite {\BG,
Corollary II.3.9}, the algebra $Z(K_n^T)$ is a Laurent polynomial ring
whose indeterminates can be chosen to be linearly independent
$H$-eigenvectors.  Thus, in calculating the center, we need only
consider an individual $H$-eigenvector, rather than a sum of such.  Once
the generators of $Z(K_n^T)$ are known, we use them to compute the
maximal ideals of the localization. Contracting these ideals will then
yield the primitive ideals of $K_n.$
 
For the case where
$n=2,$ various restrictions on the scalars
$p_1, q_1, p_2, q_2, \gamma_{1,1},\gamma_{1,2},$ and $\gamma_{2,2}$ give
rise to an assortment of primitive ideals.

The admissible sets of $K_2$ are as follows:  

$$\xalignat2  &\{x_1, y_1, \Om_1, x_2, y_2, \Om_2\}  & \{x_1, y_1,
\Om_1, y_2, \Om_2 \} \hskip1in\\ &\{x_1, y_1, \Om_1, x_2, \Om_2\}   &
\{y_1, \Om_1, x_2, y_2, \Om_2\} \hskip1in \\ &\{x_1, \Om_1, x_2, y_2,
\Om_2\} & \{y_1, \Om_1, y_2, \Om_2\} \hskip1in \\ &\{y_1, \Om_1, x_2,
\Om_2\} & \{x_1, \Om_1, x_2, \Om_2\}  \hskip1in \\ & \{x_1, \Om_1, y_2,
\Om_2\} & \{x_1, y_1, \Om_1\} \hskip1in \\ &\{y_1, \Om_1\} & \{x_1,
\Om_1\} \hskip1in \\ &\{\Om_2\} & \emptyset. \hskip1in
\endxalignat$$

\medskip
For all possible choices of the scalars $p_i,$ $q_i, \text{ and }$  
$\gamma_{i,j},$ subject only to our usual restriction that
$p_iq_i^{-1}$ is not a root of unity, each of the following ideals is 
primitive for all choices of $\alpha \in \k^{\times}$: 
$\la x_1, y_1, x_2, y_2 \ra,$  
$\la x_1, y_1, x_2 - \alpha, y_2 \ra,$ 
$\la x_1, y_1, x_2, y_2 - \alpha \ra,$ 
$\la x_1 - \alpha, y_1, x_2, y_2 \ra,$ 
and $\la x_1, y_1- \alpha, x_2, y_2 \ra.$   
This describes $\prim_{\la T \ra} K_2$ for the first five choices of the
admissible set $T$ above.  
\medskip

 The remaining primitive ideals depend more explicitly on the choices of
$p_i,$ $q_j,$ and $\gamma_{i,j}$ .   We will list some samples under the
admissible set that generates the appropriate $H$-prime ideal.

\definition{$\underline{\{y_1, \Om_1, y_2, \Om_2 \}}$}  

For $T = \{ y_1, \Om_1, y_2, \Om_2 \},$ the set $N_T = \{x_1, x_2 \}.$   
The $H$-eigenvectors in $K_n^T$ are of the form $z = \lambda\overline
x_1^a \overline x_2^b$ for some $\lambda \in k^{\times}$ and $a, b \in
\Bbb{Z}.$  We have that $z \in Z(K_n^T)$ if and only if $\overline x_1z
= z\overline x_1$ and $\overline x_2z = z \overline x_2,$ or 
$$(q_1p_2^{-1} \got)^b = 1 \text { and } (q_1p_2^{-1} \got)^a = 1.
\tag*$$
   
$\bullet$ If $\got q_1p_2^{-1}$ is not a root of unity, then $a=b = 0,$
so $z = \lambda.$  Thus,  $Z(K_n^T) = \k ,$ so $\la \overline 0 \ra$ is
the only maximal ideal of $Z(K_n^T).$  Thus, $\la y_1, y_2 \ra$ is
primitive.

$\bullet$  If $\gamma_{1,2}=1$ and $q_1p_2^{-1}$ is a root of unity 
with order $t,$ then  (*) holds if and only if $a = mt$ and $b = rt$ for
some $r, t \in \Bbb{Z}.$  Hence, $z = \lambda x_1^{mt}y_1^{rt},$ and
$Z(K_n^T) \cong \k[x_1^{\pm t}, y_1^{\pm t}].$  Then the maximal ideals
of $Z(K_n^T)$ are of the form $\la x_1^t -\alpha, x_2^t - \beta \ra$ for
$\alpha, \beta \in \k^{\times}.$  Consequently, $\la x_1^t-\alpha, y_1,
x_2^t -\beta, y_2 \ra$  is primitive for all $\alpha, \beta \in
\k^{\times}.$ 
\enddefinition

\definition{$\underline {\{y_1, \Om_1, x_2, \Om_2 \}}$}

 $\bullet$  If $p_2^{-1}\gamma_{1,2}$ is not a root of unity, then $\la y_1, x_2 \ra$ is primitive.

  $\bullet$  If $p_2^{-1}\gamma_{1,2}$ is a root of unity with order $t,$ then 
$\la x_1^t-\alpha, y_1, x_2, y_2^t-\beta \ra$ is primitive for all
$\alpha, \beta \in \k^{\times}.$

\enddefinition
 
\definition {$\underline {\{x_1, \Om_1, y_2, \Om_2 \}}$}

  $\bullet$  If $q_1\gamma_{1,2}$ is not a root of unity, then $\la x_1, y_2 \ra$ is primitive. 

  $\bullet$  If $q_1\gamma_{1,2}$ is a root of unity of order $t,$ then 
$\la x_1, y_1^t - \alpha, x_2^t - \beta, y_2 \ra$ is primitive for all 
$\alpha, \beta \in \k^{\times}.$

\enddefinition

\definition{$\underline {\{ x_1, \Om_1, x_2, \Om_2 \}}$}

$\blt$  If $\got$ is not a root of unity, then $\la x_1, x_2 \ra$ is primitive.

$\blt$  If $\got$ is a root of unity with order $t$, then $\la x_1, 
y_1^t - \alpha, x_2, y_2^t - \beta \ra$  is primitive for all $\alpha,
\beta \in \k^{\times}$.  
\enddefinition

\definition{$\underline {\{x_1, y_1, \Om_1 \}}$}

$\bullet$  If $q_2$ is not a root of unity, then $\la x_1, y_1 \ra$ is
primitive.  

$\bullet$ If $q_2$ is a root of unity with order $t,$ then for all
$\alpha, \beta \in \k^{\times},$ the ideal 
$\la x_1, y_1, x_2^t - \alpha, y_2^t - \beta \ra$ is primitive.
\enddefinition
  
\definition {$\underline {\{ \Om_2 \}}$}

$\bullet$  If $\la q_1,$ $p_2,\gamma_{1,2} \ra$ is a free abelian group
of rank 3, then $\la \Om_2 \ra$ is a primitive ideal.

 $\bullet$  If $q_1 = 1$ and $p_2 = \gamma_{1,2} $ is not a root of
unity, then $\la \Om_2, x_1 - \alpha \ra$ is primitive for all $\alpha
\in \k^{\times}.$

$\blt$  If $q_1 = 1= p_2 =\got,$ then $\la \Om_2, x_1 - \alpha, 
y_1 - \beta, x_2 - \lambda \ra$ is primitive for all $\alpha, \beta,
\lambda \in \k^{\times}.$   

\enddefinition

\definition {$\underline {\{y_1, \Om_1 \}}$}

$\blt$ If $\la q_1, p_2, \gamma_{1,2}, q_2 \ra$ is a free abelian group
of rank 4, then $\la y_1 \ra$ is primitive.

$\blt$  If $q_1 = 1=q_2$ and $\la p_2, \got \ra$ is a free abelian group
of rank 2, then $\la y_1, x_2y_2 - \alpha \ra$ is primitive for all
$\alpha \in \k^{\times}.$ 

$\blt$  If $q_1 = 1=q_2$ and $p_2 = \got $ is not a root of unity, then
for all $\alpha, \beta, \lambda \in \k^{\times}, $ the ideal
$\la y_1, x_1 - \alpha, x_2 - \beta, y_2 - \lambda \ra $ is primitive.

\enddefinition

\definition {$\underline { \{x_1, \Om_1 \} }$}

$\blt$  If $\la q_1, \got, q_2 \ra $ is a free abelian group of rank
3, then $\la x_1 \ra$ is primitive.

$\blt$  If $q_1 = 1 = q_2$ and $\got$ is not a root of unity, then $\la
x_1, x_2y_2 - \alpha \ra$ is primitive for all $\alpha \in \k^{\times}.$

$\blt$  If $q_1 = 1 = q_2 = \got,$ then $\la x_1, y_1 - \alpha, x_2 -
\beta, y_2 - \lambda \ra$ is primitive for all $\alpha, \beta, \lambda
\in \k^{\times}.$

\enddefinition

\definition {$\underline {\emptyset }$}

$\blt$  If $\la q_1, q_2, p_2, \gamma_{1,2}\ra $ is a free abelian group
of rank 4, then $\la \emptyset \ra = \la 0 \ra$ is primitive.

$\blt$  If $q_1 = p_2 = 1$ with $\la \gamma_{1,2}, q_2 \ra$ a free
abelian group of rank 2, then $\la x_1y_1 - \alpha \ra$ is primitive for
all $\alpha \in \k^{\times}.$ 

$\blt$  If $q_1=\got = q_2 = 1$ and $p_2$ is not a root of unity, then 
$\la y_1 - \alpha, \Om_2 -\beta \ra$ is primitive for all $\alpha, \beta
\in \k^{\times}.$

\enddefinition

\enddocument